\documentclass[12pt]{amsart}
\usepackage[margin=1in]{geometry}
\usepackage[cleanlook,UKenglish]{isodate}
\usepackage{ifthen}

\usepackage{subcaption}
\usepackage{enumitem}

\usepackage[backend=biber,citestyle=alphabetic,style=alphabetic,maxbibnames=9]{biblatex}

\renewbibmacro{in:}{\ifentrytype{article}{}{\printtext{\bibstring{in}\intitlepunct}}}
\DeclareFieldFormat{pages}{#1}
\DeclareFieldFormat[arxiv]{eprint}{arXiv:#1}
\DeclareFieldFormat[arxiv]{primaryClass}{}
\DeclareFieldFormat[arxiv]{secondaryClass}{}
\DeclareFieldFormat[arxiv]{archivePrefix}{}
\DeclareFieldFormat[arxiv]{title}{``#1''}
\DeclareFieldFormat{doi}{}
\DeclareFieldFormat{isbn}{}
\DeclareFieldFormat{issn}{}
\DeclareFieldFormat{url}{}
\usepackage[dvipsnames]{xcolor}
\usepackage{tikz, graphicx}
\usetikzlibrary{arrows.meta}
\usetikzlibrary{calc}
\usetikzlibrary{decorations.markings}
\usetikzlibrary{intersections}

\usepackage{amssymb,mathtools,mathabx}
\usepackage{mathrsfs}

\newcommand{\myProofSymbol}{\ensuremath{\blacksquare}}
\newcommand{\myExampleSymbol}{\ensuremath{\Box}}
\renewcommand{\qedsymbol}{\myProofSymbol}
\newcommand{\noProof}{\renewcommand{\qedsymbol}{\myProofSymbol}\qed}

\theoremstyle{plain}
\newtheorem{thm}{Theorem}[section]
\newtheorem{cor}[thm]{Corollary}
\newtheorem{lem}[thm]{Lemma}
\newtheorem{prop}[thm]{Proposition}
\newtheorem{AMSsco}[thm]{Scholium}
\newenvironment{sco}{\begin{AMSsco}}{\renewcommand{\qedsymbol}{\myProofSymbol}\qed\end{AMSsco}}
\theoremstyle{definition}
\newtheorem{AMSex}[thm]{Example}
\newenvironment{ex}{\begin{AMSex}}{\renewcommand{\qedsymbol}{\myExampleSymbol}\qed\end{AMSex}}
\theoremstyle{remark}
\newtheorem{AMSrmk}[thm]{Remark}
\newenvironment{rmk}{\begin{AMSrmk}}{\renewcommand{\qedsymbol}{\myExampleSymbol}\qed\end{AMSrmk}}

\usepackage[capitalise,nameinlink]{cleveref}
\crefname{axiom}{Axiom}{Axioms}
\crefname{figure}{Figure}{Figures}
\crefname{thm}{Theorem}{Theorems}
\crefname{cor}{Corollary}{Corollaries}
\crefname{lem}{Lemma}{Lemmas}
\crefname{prop}{Proposition}{Propositions}
\crefname{sco}{Scholium}{Scholiums}
\crefname{ex}{Example}{Examples}
\crefname{defn}{Definition}{Definitions}
\crefname{prob}{Problem}{Problems}
\crefname{rmk}{Remark}{Remarks}

\newcommand{\bb}{\mathbb}
\newcommand{\mc}{\mathcal}

\renewcommand{\check}{\widecheck}
\renewcommand{\hat}{\widehat}
\renewcommand{\setminus}{\smallsetminus}
\renewcommand{\bar}{\overline}

\DeclareMathOperator{\affineSpace}{\mathbb{E}}
\DeclareMathOperator{\R}{\mathbb{R}}
\DeclareMathOperator{\N}{\mathbb{N}_0}
\newcommand{\linearSpan}[1]{\left\langle #1 \right\rangle}

\newcommand{\set}[2]{\left\{ #1 : #2 \right\}}
\newcommand{\powerset}[2][]{\mathcal{P}_{#1}(#2)}
\newcommand{\cardinality}[1]{\left\vert{#1}\right\vert}
\DeclareMathOperator{\symmDiff}{\Delta}

\DeclareMathOperator{\id}{id}

\renewcommand{\bot}{\hat{0}}
\renewcommand{\top}{\hat{1}}

\DeclareMathOperator{\independentSetSymbol}{\mathscr{I}}

\DeclareMathOperator{\circuitSymbol}{\mathcal{C}}
\DeclareMathOperator{\flatSymbol}{\mathcal{L}}
\DeclareMathOperator{\rankSymbol}{rk}
\DeclareMathOperator{\linearClassSymbol}{\beta}
\DeclareMathOperator{\rk}{\rankSymbol}

\newcommand{\flats}[1]{\flatSymbol(#1)}

\newcommand{\circuits}[1]{\circuitSymbol(#1)}

\newcommand{\independentSets}[1]{\independentSetSymbol\left(#1\right)}
\newcommand{\fundamentalCircuit}[2]{C_{#1}(#2)}
\newcommand{\weakMap}{\leadsto}

\DeclareMathOperator{\cl}{cl}
\newcommand{\pythag}[3]{\mathcal{H}(#1; #2, #3)}

\newcommand{\systemMatrix}[3]{\operatorname{Mat}_{{#1}, {#2}}\left({#3}\right)}

\newcommand{\affineSpan}[1]{\operatorname{aff}\left(#1\right)}

\begin{document}

\title[Pythagorean Hyperplane Arrangements]{Pythagorean Hyperplane Arrangements: \\%
  Combinatorics of Gain Genericity%
}
\author{Chris Eppolito}
\address{Department of Mathematics, the University of the South, Sewanee, TN 37383, USA}
\curraddr{}
\email{christopher.eppolito@sewanee.edu}

\date{\printdate{2023-08-18}}

\begin{abstract}
  We study Pythagorean hyperplane arrangements, originally defined by Zaslavsky.
  In this first part of a series on such arrangements, we introduce a new notion of genericity for such arrangements.
  Using this notion we construct an auxiliary hyperplane arrangement whose combinatorics determines the combinatorics of all possible Pythagorean arrangements.
  We close with several applications and examples.
\end{abstract}

\maketitle

\section{Introduction}
In \cite{perpendicular_dissections_of_space--zaslavsky}, Zaslavsky initiated the study of Pythagorean arrangements.
These are hyperplane arrangements defined using two pieces of data:\ an additive real gain graph and a configuration of ``reference points'' in affine space, which are in one-to-one correspondence with vertices of the graph.
There is a hyperplane for each edge of the graph, perpendicular to the line through the points corresponding to the ends of the edge; the hyperplane's equation incorporates the gain on the edge, and determines the hyperplane's intersection with the line.

More precisely, we start with an undirected graph (with loops forbidden and parallel edges allowed) \( \Gamma \) together with a gain function \( g \) taking oriented edges of \( \Gamma \) to real numbers under the condition \( g(\bar{e}) = - g(e) \) for all oriented edges; the pair \( (\Gamma, g) \) is called an \emph{(additive real) gain graph}.
Given an additive real gain graph \( (\Gamma, g) \), a configuration \( \mc{Q} = \set{q_v}{v \in V(\Gamma)} \) of \emph{reference points}, and an edge \( e \in E(\Gamma) \) oriented from \( u \) to \( v \), we obtain a hyperplane
\begin{equation}
  h_e
  = \set{x}{d^2(x, q_u) - d^2(x, q_v) = g(e)}
  .
\end{equation}
The \emph{Pythagorean arrangement} \( \pythag{Q}{\Gamma}{g} \) associated to the triple \( (\mc{Q}, \Gamma, g) \) is the collection of all such hyperplanes.
Zaslavsky was mainly concerned with generic Pythagorean arrangements, i.e., those which have stable intersection pattern under perturbation of the reference points; that work is agnostic about what this should mean precisely, and investigates several possible notions.

Here we investigate a novel notion of genericity for Pythagorean arrangements, dubbed ``gain genericity'', which corresponds to stability of intersection patterns under shifting a hyperplane in its orthogonal direction, which is equivalent to perturbing a gain in the corresponding gain graph.
Considering this notion leads us to new results on Pythagorean arrangements in general.
Our main result is a determination of the combinatorics of every Pythagorean arrangement, \cref{result:pythagorean-arrangement-combinatorics}, which follows as a corollary of our analysis of gain-generic Pythagorean arrangements.

The main tools for our analysis are some correspondences between combinatorics and geometry.
Matroids \cite{matroid_theory--oxley, matroid_theory--welsh} and semimatroids \cite{semimatroids_and_their_tutte_polynomials--ardila, on_geometric_semilattices--wachs_walker} are the combinatorial model for hyperplane arrangements in our framework, and we further adapt results from \cite{elementary_strong_maps_between_combinatorial_geometries--dowling_kelly} to our purposes.
On the geometric side, we need only some elementary linear algebra and the related theories of real affine and projective spaces.
The combinatorics and geometry are bridged by the connection between matroids and hyperplane combinatorics, as well as by a geometric interpretation of Dilworth's truncation construction \cite{dependence_relations_in_a_semi-modular_lattice--dilworth}.

The remainder of this paper is organized as follows.
\cref{section:preliminaries} contains preliminaries and terminology used throughout.
Some attention is paid to Dilworth truncations \cite{dependence_relations_in_a_semi-modular_lattice--dilworth}, modular ideals \cite{elementary_strong_maps_between_combinatorial_geometries--dowling_kelly}, and linear classes of circuits \cite{biased_graphs_1__bias_balance_and_gains--zaslavsky, biased_graphs_2__the_three_matroids--zaslavsky}.
\cref{section:pythag} contains the meat of our results, beginning with several simple lemmas and propositions.
\cref{section:modular-ideals} further develops the relationship between linear classes of circuits and modular ideals.
\cref{section:pythagorean-combinatorics} uses these results to analyze the combinatorics of general Pythagorean arrangements, the high point being \cref{result:pythagorean-arrangement-combinatorics} and its proof from considerations of gain genericity.
Finally, \cref{section:applications} contains a variety of examples and applications of the results in \cref{section:pythag}.

This paper develops the fundamentals and tools for further generalizations of results and notions from \cite{perpendicular_dissections_of_space--zaslavsky}.
In sequels to this paper, we plan to use the results and constructions developed here to further study generic and non-generic Pythagorean arrangements from other perspectives, and with other concepts of genericity.
This will include results on extensions of a configuration-generic Pythagorean arrangement by a single reference point.

\subsection*{Acknowledgments}
This paper grew out of results in my dissertation at Binghamton University, and I'd like to thank the Department of Mathematical Sciences for their support during its writing.
Thanks are due in particular to my advisors Laura Anderson and Thomas Zaslavsky for helpful conversations, their thorough proofreading, and insightful comments.
I also thank Uly Alvarez for listening to my over-excited---and almost surely incomprehensible---rants as new discoveries occurred.

\section{Preliminaries}
\label{section:preliminaries}
Our work requires knowledge of (real) affine and projective spaces, hyperplane arrangements, matroids and their modular ideals, and (additive real) gain graphs.
This section provides a brief introduction to the relevant terminology, notation, and properties.
Readers should feel free to skip sections on topics with which they are familiar.

The set \( \N \) of natural numbers includes \( 0 \).
We let \( [n] = \{1, 2, \dots, n\} \) for each \( n \in \N \).
For a set \( S \) we let \( \binom{S}{k} = \set{T \subseteq S}{\cardinality{S} = k} \).
We might abbreviate small sets, e.g., writing \( 135 \) to denote the set \( \{1, 3, 5\} \).
The intended meaning should be clear from context.

\subsection{Affine and Projective Spaces}
By \emph{affine space} we mean a \( d \)-dimensional real affine space \( \affineSpace = \affineSpace^d \) with the usual Euclidean distance.
We often identify such spaces with \( \R^d \) for convenience.
We let \( d(x, y) \) denote the distance between points \( x, y \in \affineSpace \), and \( d^2(x, y) \) denotes the distance-squared.
An \emph{affine \( r \)-flat} of \( \affineSpace \) is a subset \( W \subseteq \affineSpace \) of points so that the set \( D(W) = \set{q - p}{p, q \in W} \) of \emph{directions in \( W \)} forms a real vector space of dimension \( r \).
Fixing a base point \( p \in W \) specifies an origin, and we may express \( D(W) = \set{q-p}{q \in W} \).
A \emph{hyperplane} is an affine \( (d-1) \)-flat, a \emph{line} is an affine \( 1 \)-flat, and a \emph{point} is an affine \( 0 \)-flat.
The empty set is the unique affine \( (-1) \)-flat.
Affine flats are \emph{parallel} when they are equal or have empty intersection.

The \emph{affine span} of \( S \subseteq \affineSpace \) is \( \affineSpan{S} \), the intersection of all affine subspaces of \( \affineSpace \) containing \( S \).
A point \( p \) is \emph{affinely dependent} on \( S \subseteq \affineSpace \) when \( p \in \affineSpan{S} \), i.e., there are a finite \( T \subseteq S \) and coefficients \( c_t \in \R \) with \( \sum_{t \in T} c_t = 1 \) and \( \sum_{t \in T} c_t(t - p) = \vec{0} \).
Such an \emph{affine dependence relation} is also written \( p = \sum_{t \in T}c_t t \).

Every affine space \( \affineSpace \) has a \emph{projective completion} \( \hat{\affineSpace} \) obtained by appending a ``hyperplane at infinity''.
Formally, a \emph{point at infinity} is an equivalence class of lines in \( \affineSpace \) modulo the parallelism relation, i.e., \( l_1 \parallel l_2 \) when there are \( p_1, q_1 \in l_1 \) and \( p_2, q_2 \in l_2 \) with \( p_i \neq q_i \) and \( q_1 - p_1 = q_2 - p_2 \) in \( D(\affineSpace) \).
The hyperplane at infinity is the set \( H_\infty = \set{l/_\parallel}{l \text{ is a line in } \affineSpace} \), where \( l/_\parallel \) denotes the equivalence class of \( l \).
The projective completion is a projective space; indeed, each affine \( r \)-flat \( A \subseteq \affineSpace \) yields a \emph{projective \( r \)-flat} \( \hat{A} = A \cup \set{l/_\parallel}{l \text{ is a line in } A} \).
We often use the following fact.
\begin{lem}
  If \( P \) is a projective space and \( H \) is a hyperplane therein, then \( P \setminus H \) is an affine space.
  Furthermore, \( P \setminus H \) is isomorphic to \( P \setminus H' \) for all \( H \) and \( H' \).
  \noProof
\end{lem}
This lemma allows us to treat any hyperplane \( H \) in a projective space \( P \) as the hyperplane at infinity in the projective completion of the resulting affine space \( P \setminus H \).
Moreover, because we work in real affine (resp.\ projective) spaces, given any finite collection of points it is always possible to find an affine (resp.\ projective) hyperplane containing none of the given points.
This allows us to treat finitely many points in the hyperplane at infinity as points in an affine space, a trick we use below.

\subsection{Hyperplane Arrangements}
A \emph{hyperplane arrangement} is a finite multiset \( \mc{A} \) of hyperplanes in affine space \( \affineSpace \).
A multisubset of \( \mc{A} \) is called a \emph{subarrangement} of \( \mc{A} \).
An arrangement is \emph{essential} when there is a subset \( S \subseteq \mc{A} \) so that \( \bigcap S \) is a point, and it is \emph{central} when \( \bigcap \mc{A} \) is nonempty.
We sometimes talk about the ``linear system associated to a subarrangement'', by which we mean the linear system of the hyperplane equations in that subarrangement
The \emph{intersection semilattice} of \( \mc{A} \) is the set of all nonempty intersections of multisubsets of \( \mc{A} \), ordered by reverse inclusion; these intersections are the \emph{flats} of the arrangement.
Two hyperplane arrangements have the \emph{same combinatorial type} or \emph{same intersection pattern} when they have isomorphic intersection semilattices.\footnote{As distinct hyperplanes may coincide as sets, we require an isomorphism of the labeled intersection semilattices, i.e., one induced by a bijection of label sets; we often suppress this technicality.}
The intersection semilattice is ranked by co-dimension.

\subsection{Gain Graphs}
A \emph{graph} \( \Gamma \) has a set of \emph{vertices} and a set of edges \emph{edges}, denoted by \( V(\Gamma) \) and \( E(\Gamma) \) respectively.
Each edge has two \emph{end vertices}, or \emph{ends}.
A \emph{loop} is an edge with both ends at the same vertex, and edges \( a, b \) are \emph{parallel} when they have the same multiset of ends; a graph is \emph{simple} when it is has neither loops nor parallel edges.
A vertex is \emph{isolated} when it is not an end of any edge.
We write \( \Delta \subseteq \Gamma \) to denote that \( \Delta \) is a \emph{subgraph} of \( \Gamma \), i.e., that \( \Delta \) is a graph satisfying \( V(\Delta) \subseteq V(\Gamma) \) and \( E(\Delta) \subseteq E(\Gamma) \).
For \( S \subseteq V(\Gamma) \), let \( \Gamma[S] \) denote the \emph{induced subgraph} on vertex set \( S \), i.e., the maximal subgraph of \( \Gamma \) with \( V(\Gamma[S]) = S \).
We often denote \( \Gamma[V(\Gamma) \setminus \{v\}] \) by \( \Gamma \setminus v \) for brevity.
We sometimes conflate a subgraph and its edge set; our work does not change in the presence of isolated vertices, so this is harmless.

A \emph{walk} in \( \Gamma \) is a sequence \( (v_0, e_1, v_1, \dots, v_{k-1}, e_k, v_k) \) of vertices \( v_i \) and edges \( e_i \) with ends \( v_{i-1} \) and \( v_i \).
Each edge has two possible \emph{orientations} corresponding to the directions we might traverse the edge in a walk.
A \emph{closed walk} is a walk with the same initial and terminal vertex.
A \emph{circle} is a (sub)graph which can be traced by an internally non-repeating closed walk; we sometimes treat a circle as one of its corresponding closed walks.
A \emph{forest} is a graph without circles.
The \emph{rank} of a subgraph is the maximum size of a forest contained therein, measured by number of edges.

Let \( G \) be a (multiplicative) group.
A \emph{\( G \)-gain graph} is a graph \( \Gamma \) with a \emph{gain function} \( g \) mapping oriented edges of \( \Gamma \) to elements of \( G \), obeying the identity \( g(\bar{e}) = g(e)^{-1} \) where \( \bar{e} \) denotes \( e \) with the opposite orientation.
The gain function canonically extends to a function taking walks in \( \Gamma \) to elements of \( G \) by taking the product of the gains on the edges in the order visited.
A walk \( W \) is \emph{neutral} when \( g(W) = 1_G \); a graph is \emph{balanced} when all closed walks in the subgraph are neutral.
Equivalently, a graph is balanced when all circle subgraphs contained therein are balanced; see \cite[Lemma 5.3]{biased_graphs_1__bias_balance_and_gains--zaslavsky}.

The \emph{balance-closure} of a subset \( S \subseteq E(\Gamma) \) is
\[
  \cl S
  = S \cup \set{e \in E(\Gamma)}
  {\text{there is a balanced circle } C \text{ with } C \subseteq S \cup \{e\}}
  .
\]
\begin{lem}
  The balance-closure of a balanced set \( S \) is the largest balanced set containing \( S \) and having the same rank in the cycle matroid of the graph.
\end{lem}
\begin{proof}
  Suppose \( T \supseteq S \) is balanced, \( \rk T = \rk S \), and \( e \in T \supseteq S \).
  As their ranks are the same, there is a maximal forest of \( S \) which is a maximal forest of \( T \).
  Thus the ends of \( e \) are connected by a path \( P \) in \( S \), which determines the gain on \( e \) by balance.
  Now \( C = P \cup \{e\} \) is a circle with \( e \in C \subseteq S \cup \{e\} \), which yields \( e \in \cl S \).
\end{proof}
All forests are balanced, so their balance-closure is also balanced.

A \emph{theta graph} is a union of three internally disjoint paths between a pair of distinct vertices.
A subset \( \linearClassSymbol \) of the circles of a graph is a \emph{linear class of circles} when for all \( C_1, C_2 \in \linearClassSymbol \) with \( C_1 \cup C_2 \) a theta graph, the circle \( C_1 \symmDiff C_2 \) is also in \( \linearClassSymbol \).

\subsection{Matroids}
Let \( E \) be a finite set.
A collection \( \circuitSymbol \) of nonempty, pairwise incomparable subsets of \( E \) is the set of \emph{circuits} of a matroid on \( E \) when it satisfies the \emph{(weak) circuit elimination} axiom, i.e.,
\begin{enumerate}
\item[(CW)]\label[axiom]{circuit:weak-elimination} %
  For all circuits \( X \) and \( Y \) and all \( e \in X \cap Y \) there is a circuit \( Z \) such that \( Z \subseteq (X \cup Y) \setminus \{e\} \).
\end{enumerate}
For a collection of pairwise incomparable and nonempty sets, this is equivalent to the following \emph{a priori} stronger condition---the \emph{strong circuit elimination} axiom.
\begin{enumerate}
\item[(CS)]\label[axiom]{circuit:strong-elimination} %
  For all circuits \( X \) and \( Y \), all \( x \in X \cap Y \), and all \( y \in Y \setminus X \) there is a circuit \( Z \) such that \( y \in Z \subseteq (X \cup Y) \setminus \{e\} \).
\end{enumerate}
A set is \emph{dependent} when it contains a circuit, and \emph{independent} otherwise.
A \emph{basis} is a maximal independent set.
Given a basis \( B \) and an \( e \in E \setminus B \), there is a unique \emph{fundamental circuit} \( \fundamentalCircuit{B}{e} \) of \( e \) with respect to \( B \) satisfying \( e \in \fundamentalCircuit{B}{e} \subseteq B \cup \{e\} \).

The \emph{rank} of a subset $S \subseteq E$ is the size of a maximal independent set contained therein.
The \emph{closure} of $S$ is the maximal subset $\sigma(S) \subseteq E$ containing $S$ and for which $\rk(\sigma(S)) = \rk(S)$; we have
\[
  \sigma(S) = S \cup \set{e \in E}{C \subseteq S \cup e \text{ for some } C \in \mc{C}}
  .
\]
The dual of a matroid \( M \) is the matroid \( M^* \) whose bases are the complements of the bases of \( M \).

A \emph{loop} of $M$ is an element $e \in E$ for which $\{e\} \in \mc{C}$.
Two elements \( x, y \) of a matroid are \emph{parallel} when they are either equal, both loops, or belong to a common circuit of size two.
Parallelism is an equivalence relation on the ground set of a matroid.
A matroid is \emph{simple} when it has neither loops nor nontrivial parallelism.

\begin{figure}
  \centering
  \pgfmathsetmacro{\loX}{0}
  \pgfmathsetmacro{\hiX}{6}
  \pgfmathsetmacro{\loY}{0}
  \pgfmathsetmacro{\midY}{3}
  \pgfmathsetmacro{\hiY}{4}
    \begin{tikzpicture}[scale=.7, transform shape]
      \path
      (\loX,\loY) coordinate (4)
      (\hiX,\loY) coordinate (8)
      (4) -- (8) coordinate[pos=.5] (3)
      (\loX,\midY) coordinate (2)
      (\hiX,\hiY) coordinate (6)
      (2) -- (6) coordinate[pos=.5] (7)
      (intersection of 2--3 and 4--7) coordinate (9)
      (intersection of 6--3 and 7--8) coordinate (1)
      (intersection of 2--8 and 6--4) coordinate (5)
      ;
      \draw[every node/.style={circle, inner sep=1.5pt, fill=black}]
      (1) node[label={ 90:$1$}] {}
      (2) node[label={ 90:$2$}] {}
      (3) node[label={270:$3$}] {}
      (4) node[label={270:$4$}] {}
      (5) node[label={ 90:$5$}] {}
      (6) node[label={ 90:$6$}] {}
      (7) node[label={ 90:$7$}] {}
      (8) node[label={270:$8$}] {}
      (9) node[label={90:$9$}] {}
      ;
      \foreach \i/\j in {2/3,2/8,7/4,7/8,6/4,6/3,2/6,4/8}{
        \draw[black] ($ (\i)!-.25!(\j) $) -- ($ (\i)!1.25!(\j) $);
      }
    \end{tikzpicture}
  \caption{%
    \label{figure:pappos-related-matroids}%
    Non-Pappos matroid.
  }
\end{figure}
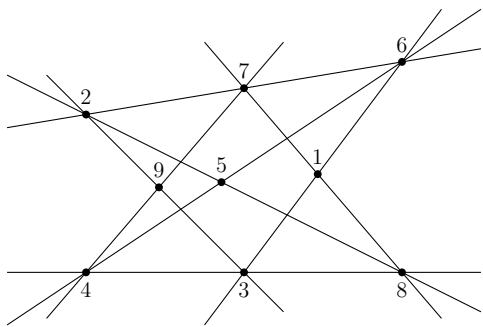

We sometimes represent a matroid of rank \( \leq 3 \) with an \emph{affine diagram}.
In such a diagram, e.g., \cref{figure:pappos-related-matroids}, the elements of the matroid are represented as points, and circuits of rank \( 2 \) are represented as collinearities (though we may draw lines which are not straight to achieve this effect if necessary).
Note that circuits of rank \( 3 \) are not depicted; these are exactly the subsets of size \( 4 \) containing no collinear triple.

\begin{ex}
Every graph \( \Gamma \) has a \emph{cycle matroid} with ground set \( E(\Gamma) \) and whose circuits are the circles of \( \Gamma \); such matroids are \emph{graphic}.
Every finite collection \( E \) of vectors in a vector space has a \emph{vectorial matroid} with ground set \( E \) and whose independent sets are the linearly independent subsets of \( E \).
\end{ex}

A \emph{weak map} \( f \colon M \weakMap N \) of matroids is a function \( f \colon E(M) \to E(N) \) so that the image of every circuit of \( M \) is dependent in \( N \).
If \( E(M) = E = E(N) \) and \( f = \id_E \), then \( f \) is a weak map if and only if every circuit of \( N \) contains a circuit of \( M \).
Thus, weak maps are the matroid analogue of ``specialization of position'' (see \cref{figure:hyperplane-weak-map}).
\begin{figure}
  \centering
  \begin{tikzpicture}[>=Stealth]
    \pgfmathsetmacro{\myRadius}{2}
    \pgfmathsetmacro{\vecRadius}{1}
    \pgfmathsetmacro{\shiftParam}{\myRadius+2}
    \pgfmathsetmacro{\angleOne}{0}
    \pgfmathsetmacro{\angleTwo}{60}
    \pgfmathsetmacro{\angleThree}{210}
    \begin{scope}
      \foreach \ang/\lab/\col in {\angleOne/1/black,\angleTwo/2/blue,\angleThree/3/red}{
        \pgfmathsetmacro{\lAngle}{\ang+90}
        \pgfmathsetmacro{\rAngle}{\ang+270}
        \coordinate (\lab) at (\ang:\vecRadius);
        \draw[\col]
        (\lAngle:\myRadius) edge[thick] node[pos=0,label={\lAngle:$h_{\lab}$}] {} (\rAngle:\myRadius)
        (0,0) edge[->] node[pos=1, label={\ang:$v_{\lab}$}] {} (\lab)
        ;
      }
    \end{scope}
    \draw (\shiftParam,0) node {$\weakMap$};
    \begin{scope}[shift={(2*\shiftParam,0)}]
      \pgfmathsetmacro{\angleThree}{\angleTwo+180}
      \foreach \ang/\lab/\col in {\angleOne/1/black,\angleTwo/2/purple,\angleThree/3/purple}{
        \pgfmathsetmacro{\lAngle}{\ang+90}
        \pgfmathsetmacro{\rAngle}{\ang+270}
        \coordinate (\lab) at (\ang:\vecRadius);
        \draw[\col]
        (\lAngle:\myRadius) edge[thick] node[pos=0,label={\lAngle:$h_{\lab}$}] {} (\rAngle:\myRadius)
        (0,0) edge[->] node[pos=1, label={\ang:$v_{\lab}$}] {} (\lab)
        ;
      }
    \end{scope}
  \end{tikzpicture}
  \caption{%
    \label{figure:hyperplane-weak-map}%
    Weak map \( f = \id_E \colon M \weakMap N \) of matroids.
  }
\end{figure}

An \emph{elementary quotient} of \( M \) is a matroid \( N \) for which \( E(M) = E = E(N) \), \( \rk N = \rk M - 1 \), and every flat of \( N \) is a flat of \( M \).
In this case, every flat of \( M^* \) is a flat of \( N^* \) and \( \rk N^* = \rk M^* + 1 \); we call \( N^* \) an \emph{elementary lift} of \( M^* \).

A pair of (distinct) sets \( S, T \subseteq E \) is a \emph{modular pair} in \( M \) when
\[
  \rk(S \cap T) + \rk(S \cup T)
  = \rk S + \rk T
  .
\]
We shall be particularly interested in modular pairs of circuits.
\begin{lem}[{\cite[Lemma 3.7]{matroids_over_partial_hyperstructures--baker_bowler}}]
  \label{result:equivalent-descriptions-of-modular-pairs}
  For distinct \( X, Y \in \circuits{M} \), the following are equivalent.
  \begin{enumerate}
  \item\label[part]{mod-pair} %
    The pair \( X, Y \) is a modular pair.
  \item\label[part]{usual} %
    \( \rk(X \cap Y) + \rk(X \cup Y) = \rk(X) + \rk(Y) \).
  \item\label[part]{minus-2} %
    \( \rk(X \cup Y) = \cardinality{X \cup Y} - 2 \).
  \item\label[part]{basis} %
    There are a basis \( B \) and \( x \in X \), \( y \in Y \) such that \( X = C_B(x) \) and \( Y = C_B(y) \).
  \item\label[part]{minimal} %
    For all distinct \( U, V \in \circuits{M} \), if \( U \cup V \subseteq X \cup Y \), then \( U \cup V = X \cup Y \).
    \noProof
  \end{enumerate}
\end{lem}

Applying \cref{result:equivalent-descriptions-of-modular-pairs} and some elementary graph theory yields the following.
\begin{cor}
  Let \( \Gamma \) be a graph with given circles \( X \) and \( Y \).
  In the cycle matroid \( M[\Gamma] \), a pair \( X, Y \) of circles forms a modular pair in the cycle matroid \( M[\Gamma] \) precisely when one of the following holds.
  \begin{enumerate}
  \item %
    The subgraph \( X \cup Y \) is disconnected.
  \item %
    The circles \( X \) and \( Y \) share a unique common vertex in \( \Gamma \).
  \item %
    The subgraph \( X \cup Y \) is a \( \theta \)-graph in \( \Gamma \).
    \noProof
  \end{enumerate}
\end{cor}

A \emph{modular ideal} of a matroid \( M \) on \( E \) is a nonempty order ideal \( I \) in \( (\powerset{E}, \subseteq) \) such that
\begin{enumerate}[label={(MI\arabic*)}]
\item %
  For all non-loops \( e \in E \) we have \( \{e\} \in \mc{I} \).
  \hfill
  (Non-degeneracy)
\item %
  For all modular pairs \( A, B \in \mc{I} \) we have \( A \cup B \in \mc{I} \).
  \hfill
  (Modular Extension)
\end{enumerate}

\begin{rmk}
  Note that \cite{elementary_strong_maps_between_combinatorial_geometries--dowling_kelly} requires that a modular ideal be a \emph{proper} order ideal; for the applications they have in mind, this is sufficient.
  We will not include this restriction, as a modular ideal will model the collection of central subsets of a hyperplane arrangement, and the full arrangement can be central.
  Their theorems are trivially altered by this change; in effect, we append a top element \( \top = \powerset{E}\) to their semilattice of modular ideals.
  Making corresponding changes elsewhere in their definitions preserves their results.
\end{rmk}

\begin{ex}
  \label{example:unrealizable-lift-matroid}
  The non-Pappos matroid from \cref{figure:pappos-related-matroids} is an elementary lift of the rank two uniform matroid on nine elements, \( U_{2, 7} \).
  Let
  \[
    \mc{X} = \{136, 178, 239, 258, 267, 348, 456, 479\}
  \]
  and \( \mc{I} = \binom{[9]}{2} \cup \mc{X} \).
  Then \( \mc{I} \) is a modular ideal of \( U_{2, 7} \); indeed, if \( A, B \in \mc{I} \) form a modular pair, then they are either comparable or both independent as \( \cardinality{A \cap B} \leq 1 \) for all \( A, B \in \mc{X} \).
  The corresponding elementary lift derived from this modular ideal via \cite[Proposition 6.4]{elementary_strong_maps_between_combinatorial_geometries--dowling_kelly} is the non-Pappos matroid.
  The elements of \( \mc{X} \) are the lines illustrated in \cref{figure:pappos-related-matroids}.

  In addition, this example shows that an elementary lift of a matroid which is realizable by a vector arrangement need not be realizable.
  Indeed, \( U_{2, 7} \) is realizable using any seven distinct vectors of the form \( (t, 1) \) in \( \R^2 \).
  On the other hand, by Pappos' Theorem in projective geometry, the non-Pappos matroid is not realizable.
\end{ex}

Note that every independent subset of \( M \) belongs to every modular ideal of \( M \).
Indeed, given an independent set \( S = \set{s_i}{1 \leq i \leq n} \) we let \( S_k = \set{s_i}{1 \leq i \leq k} \) for \( 0 \leq i \leq n \).
Thus \( S_0 = \emptyset \) and \( S_1 = \{s_1\} \) are both elements of \( \mc{I} \) as \( \mc{I} \) is a non-degenerate ideal.
As \( S_k = S_{k-1} \cup \{s_k\} \) and \( \rk_M(S_k) = k = (k-1) + 1 = \rk_M(S_{k-1}) + \rk_M\{s_k\} \) we see that \( S_k = S_{k-1} \cup \{s_k\} \in \mc{I} \) by the modular extension axiom.

\begin{prop}
  [{\cite[Corollary 6.2.2 with \( \hat{1} \)]{elementary_strong_maps_between_combinatorial_geometries--dowling_kelly}}]
  The set of modular ideals of a matroid \( M \) is a lattice under inclusion; the meet operation is intersection.
  The zero of this lattice is the ideal of independent sets of \( M \).
  \noProof
\end{prop}

Moreover, the proper modular ideals of \( M \) are in one-to-one correspondence with the so-called \emph{elementary lifts} of \( M \)---these are called \emph{elementary preimages} in \cite{elementary_strong_maps_between_combinatorial_geometries--dowling_kelly}.
\begin{prop}
  [{\cite[Proposition 6.4]{elementary_strong_maps_between_combinatorial_geometries--dowling_kelly}}]
  Let \( \mc{I} \) be a modular ideal of \( M \).
  The function
  \[
  \rk_{\mc{I}}(S)
  =
  \begin{cases}
    \rk_M(S)    & \text{if } S \in \mc{I}, \\
    \rk_M(S) + 1 & \text{otherwise}.
  \end{cases}
  \]
  is the rank function of a matroid.
  \noProof
\end{prop}
The matroid with rank function \( \rk_{\mc{I}} \) is the \emph{elementary \( \mc{I} \)-lift} of \( M \).
Allowing the powerset to be a modular ideal effectively allows construction of the ``trivial lift'' of \( M \), i.e., \( M \) itself.

\begin{prop}
  [{\cite[Proposition 6.5]{elementary_strong_maps_between_combinatorial_geometries--dowling_kelly}}]
  The lattice of modular ideals is isomorphic to the lattice of elementary lifts of \( M \) ordered by weak map.
  \noProof
\end{prop}

A hyperplane arrangement determines a matroid when it is central; otherwise, it determines a \emph{semimatroid}.
As with matroids, semimatroids admit a variety of equivalent descriptions.
Perhaps the easiest description is via \cite{on_geometric_semilattices--wachs_walker}.
A set \( \mc{F} \) of subsets of \( E \) is the collection of \emph{flats} of a semimatroid on \( E \) when there is a matroid \( M \) on \( E \cup \{z_0\} \) for which \( \mc{F} \) is the set of flats of \( M \) not containing \( z_0 \).

A semimatroid can be thought of as a ``cross-section'' of a matroid.
The geometric intuition for this description comes from hyperplane arrangements.
Let \( \mc{A} \) be a linear arrangement in \( \R^{d+1} \) (i.e., the origin is in every hyperplane), and assume there is exactly one hyperplane \( z_0 \) in \( \mc{A} \) with defining equation \( x_{d+1} = 0 \).
Now \( \bb{A} = \set{x \in \R^{d+1}}{x_{d+1} = 1} \) is an affine hyperplane, thus an affine space of dimension \( d \).
Moreover, every element of \( \mc{A} \setminus \{z_0\} \) intersects \( \bb{A} \) in an affine hyperplane thereof, so we obtain an arrangement \( \check{\mc{A}} = \set{h \cap \bb{A}}{z_0 \neq h \in \mc{A}} \) in \( \bb{A} \).
A set \( C \subseteq \check{\mc{A}} \) is central if and only if \( \bigcap C \setminus z_0 \neq \emptyset \).
In particular, the flats of \( \check{\mc{A}} \) are in bijection with the flats of \( \mc{A} \) not contained in \( z_0 \), so flats with \( z_0 \notin F \subseteq \mc{A} \) in the matroid of \( \mc{A} \).

This geometric intuition also works in reverse; given an affine arrangement in \( \affineSpace \), use any isomorphism to identify \( \affineSpace \) with \( \bb{A} \) and define \( \hat{A} = \{z_0\} \cup \set{\linearSpan{h}}{h \in \mc{A}} \).
Then \( \hat{\mc{A}} \) is a linear arrangement with \( \mc{A} = \check{\hat{\mc{A}}} \).

Following \cite{semimatroids_and_their_tutte_polynomials--ardila}, we define semimatroids via a collection \( \mc{S} \) of \emph{central sets} and a \emph{rank function} \( \rk \).
A pair \( (\mc{S}, \rk \colon \mc{S} \to \N) \) determines a \emph{semimatroid} on \( E \) when...
\begin{enumerate}[label={(S\arabic*)}, ref={S\arabic*}, start=0]
\item %
  \(\mc{S}\) is a nonempty, hereditary set of subsets of \( E \).
\item %
  \( \rk \) is subcardinal.
\item %
  \( \rk \) is isotone.
\item %
  If \( X \), \( Y \), and \( X \cup Y \) are central, then \( \rk(X \cap Y) + \rk(X \cup Y) \leq \rk X + \rk Y \).
\item %
  If \( X \) and \( Y \) are central and \( \rk(X \cap Y) = \rk X \), then \( X \cup Y \) is central.
\item %
  If \( X \) and \( Y \) are central and \( \rk X < \rk Y \), then \( X \cup \{y\} \) is central for some \( y \in Y \setminus X \).
\end{enumerate}
This definition is not the original definition of a semimatroid.
However, it is equivalent and matches our intuition from hyperplane arrangements most closely among the known descriptions, as evidenced by the next example.

\begin{ex}
  The central subsets of a hyperplane arrangement are the central sets of a semimatroid, ranked by codimension of the intersection.

  The balanced edge subsets of a gain graph form the central sets of a semimatroid.
  The rank of a balanced set is its rank in the cycle matroid of the underlying graph.
  This is implicit in \cite[Theorem 2.1]{biased_graphs_2__the_three_matroids--zaslavsky}.
\end{ex}

\begin{prop}
  [{\cite[Propositions 3.4--3.5]{semimatroids_and_their_tutte_polynomials--ardila}}]
  \label{lemma:semimatroid-iff-modular-ideal}
  Let \( \mc{S} \) be an order ideal of subsets of \( E \), and let \( \rk \colon \mc{S} \to \N \).
  Then \( \mc{S} \) is the set of central subsets of a semimatroid on \( E \) ranked by \( \rk \) if and only if \( \mc{S} \) is a modular ideal of a matroid \( M \) with rank function \( \hat{\rk} \) with \( \hat{\rk}|_{\mc{S}} = \rk \).
  \noProof
\end{prop}

\subsection{Point Configurations}
A \emph{point configuration} in \( \affineSpace \) is a finite set \( \mc{Q} \subseteq \affineSpace \) of distinct points therein.
A configuration is in \emph{affine (general) position} when every subset of \( 0 \leq n \leq d+1 \) points spans an \( n-1 \) flat.
We may project the points of \( \mc{Q} \) to infinity to obtain a point arrangement \( \mc{Q}_\infty = \set{\affineSpan{p, q}/_\parallel}{p, q \in \mc{Q} \text{ distinct}} \) in the hyperplane at infinity.
The configuration \( \mc{Q} \) is in \emph{ideal (general) position} when for every oriented forest \( F \) with \( d \) edges in the complete graph on vertices \( V = \mc{Q} \), the set
\[
  B_F = \set{t - s}{e \in E(F) \text{ is oriented from } s \text{ to } t}
\]
is linearly independent in the direction space of \( \affineSpace \).
Ideal position is stronger than affine position---a square in \( \affineSpace^2 \) has vertices in affine position but not ideal position.

Let \( \mc{Q} \) be a finite point set in affine \( d \)-space, and pick an arbitrary orientation of the complete graph on vertex set \( \mc{Q} \).
The \emph{matroid of \( \mc{Q} \) at infinity} is the matroid \( M_\infty(\mc{Q}) \) of the vector arrangement
\[
  E(M_\infty(\mc{Q}))
  = \set{q_j - q_i}{ij \text{ is oriented from } i \text{ to } j}
  .
\]
Note that \( S \subseteq E(M_\infty(\mc{Q})) \) is linearly independent if and only if \( \dim (\bigcap_{s \in S} s^\perp) = d - \cardinality{S} \).
Moreover, this condition is independent of the ordering on \( \mc{Q} \); reordering can only introduce some minus signs, and cannot change the support of a dependency.

\section{Pythagorean Arrangements}
\label{section:pythag}

\begin{figure}
  \centering
  \newcommand{\gAB}{0}
  \newcommand{\gAC}{-3}
  \newcommand{\gAD}{-2}
  \newcommand{\gBD}{-6}
  \newcommand{\gCD}{1}

  \begin{subfigure}{.65\textwidth}
    \centering
    \begin{tikzpicture}[scale=0.75]
      \pgfmathsetmacro{\xMin}{-3}
      \pgfmathsetmacro{\xMax}{10}
      \pgfmathsetmacro{\yMin}{-2}
      \pgfmathsetmacro{\yMax}{10}
      \clip (\xMin,\yMin) rectangle (\xMax,\yMax);
      \draw[every node/.style={circle, fill=black, inner sep=1pt}]
      ( 0,  0) node[label={180:$q_1$}] (A) {}
      ( 2,  3) node[label={270:$q_2$}] (B) {}
      (-1,  7) node[label={  0:$q_3$}] (C) {}
      ( 9,  2) node[label={180:$q_4$}] (D) {}
      ;
      \begin{scope}
        [thick,
        domain=\xMin:\xMax,
        plotlabel/.style args={at #1 #2 with #3}{%
          decoration={markings, mark={at position #1 with \node[#2] {#3};}},
          postaction={decorate}
        }
        ]
        \draw[plotlabel=at 0 above right with {$h_{a=12}$}] plot(\x,{(\gAB + 13 - 4*\x)/6});
        \draw[plotlabel=at 0 below right with {$h_{b=13}$}] plot(\x,{(\gAC + 50 - -2*\x)/14});
        \draw[plotlabel=at .425 left with {$h_{c=14}$}] plot(\x,{(\gAD + 85 - 18*\x)/4});
        \draw[plotlabel=at .697 left with {$h_{d=24}$}] plot(\x,{(\gBD + 72 - 14*\x)/-2});
        \draw[plotlabel=at .73 right with {$h_{e=34}$}] plot(\x,{(\gCD + 35 - 20*\x)/-10});
      \end{scope}
    \end{tikzpicture}
    \subcaption{Arrangement.}
  \end{subfigure}
  \vspace{2ex}

  \begin{subfigure}[b]{.25\textwidth}
    \centering
    \begin{tikzpicture}
      \draw[every node/.style={circle, fill, inner sep=1.5pt}]
      (2,0) node[label={270:$1$}] (A) {}
      (2,2) node[label={135:$2$}] (B) {}
      (4,0) node[label={315:$3$}] (C) {}
      (4,2) node[label={ 45:$4$}] (D) {}
      ;
      \draw[
      >=Stealth,
      ->/.style={decoration={markings, mark={at position .5 with {\arrow{>}}}}, postaction={decorate}},
      ]
      (A) edge node[left]            {$\gAB$} (B)
      (A) edge[->] node[below]       {$\gAC$} (C)
      (A) edge[->] node[below right] {$\gAD$} (D)
      (B) edge[->] node[above]       {$\gBD$} (D)
      (C) edge[->] node[right]       {$\gCD$} (D)
      ;
    \end{tikzpicture}
    \subcaption{Gain graph.}
  \end{subfigure}
  \begin{subfigure}[b]{.7\textwidth}
    \centering
    \begin{tikzpicture}[xscale=.5, yscale=1.5, every node/.style={circle, inner sep=0pt}]
      \draw
      (-8,0) node (a)   {$a$}
      (-4,0) node (b)   {$b$}
      ( 0,0) node (c)   {$c$}
      ( 4,0) node (d)   {$d$}
      ( 8,0) node (e)   {$e$}
      (-8,2) node (ab)  {$ab$}
      (-6,2) node (ace) {$ace$}
      (-4,2) node (ad)  {$ad$}
      (-2,2) node (bc)  {$bc$}
      ( 0,2) node (bd)  {$bd$}
      ( 2,2) node (be)  {$be$}
      ( 4,2) node (cd)  {$cd$}
      ( 6,2) node (ce)  {$ce$}
      ( 8,2) node (de)  {$de$}
      ;
      \draw
      (ab.south) edge (a.north) edge (b.north)
      (ace.south) edge (a.north) edge (c.north) edge (e.north)
      (ad.south) edge (a.north) edge (d.north)
      (bc.south) edge (b.north) edge (c.north)
      (bd.south) edge (b.north) edge (d.north)
      (be.south) edge (b.north) edge (e.north)
      (cd.south) edge (c.north) edge (d.north)
      (ce.south) edge (c.north) edge (e.north)
      (de.south) edge (d.north) edge (e.north)
      ;
    \end{tikzpicture}
    \subcaption{intersection semilattice (with \( \bot \) omitted).}
  \end{subfigure}

  \caption{%
    \label{figure:pythagorean-arrangement-and-intersection-poset}%
    Pythagorean arrangement and its intersection semilattice.
  }
\end{figure}
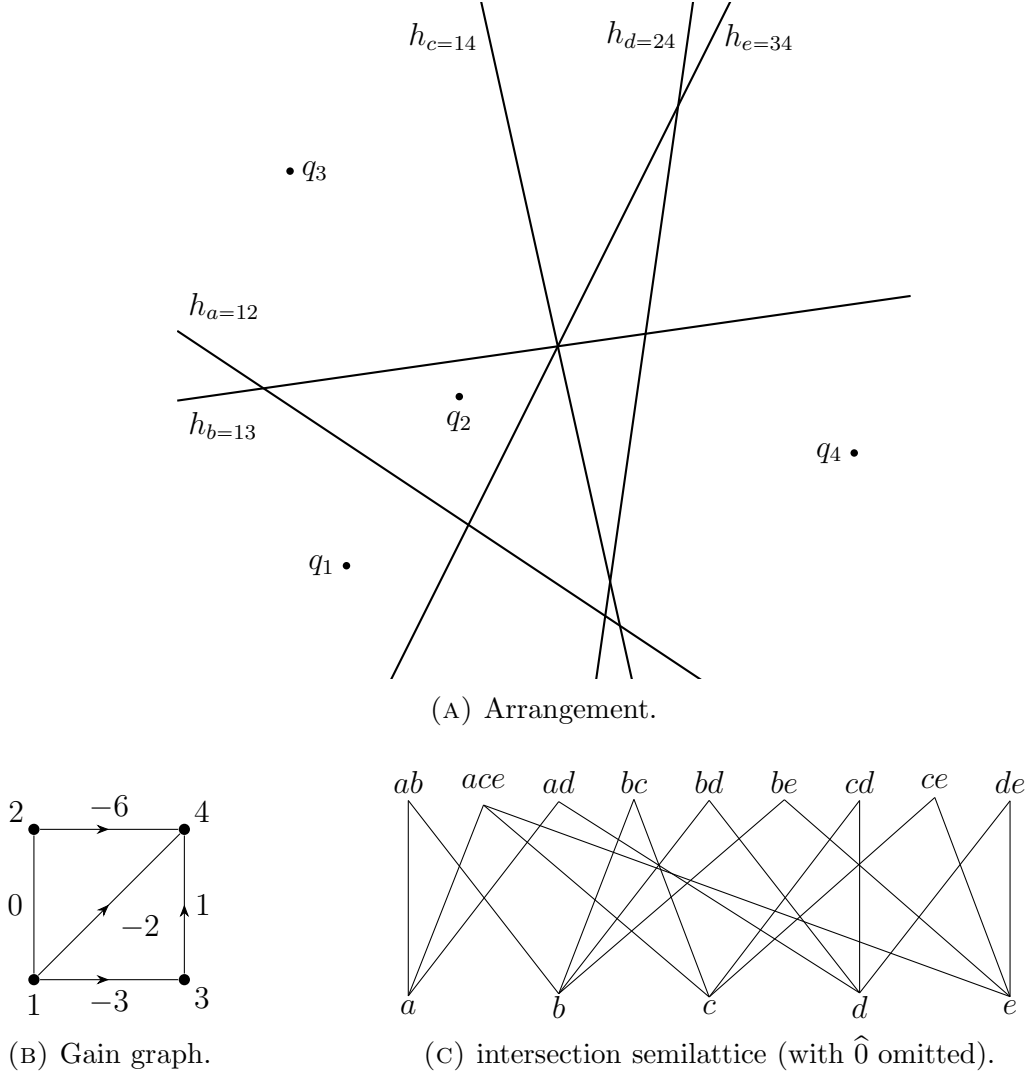

Given an additive real gain graph \( (\Gamma, g) \), a configuration \( \mc{Q} = \set{q_v}{v \in V(\Gamma)} \) of \emph{reference points}, and an edge \( e \in E(\Gamma) \) oriented from \( u \) to \( v \), we obtain a hyperplane
\begin{equation}
  h_e
  = \set{x}{d^2(x, q_u) - d^2(x, q_v) = g(e)}
  .
\end{equation}
The \emph{Pythagorean arrangement} \( \pythag{Q}{\Gamma}{g} \) associated to the triple \( (\mc{Q}, \Gamma, g) \) is the collection of all such hyperplanes.
\cref{figure:pythagorean-arrangement-and-intersection-poset} illustrates a Pythagorean arrangement.

\begin{rmk}
  The requirement that \( \mc{Q} \) have no repeated points is a simplifying assumption.
  All of our results admit trivial extensions to cases involving repeated points in one of several ways.
  Either
  \begin{enumerate}
  \item we require that the association of points with vertices is a proper vertex-coloring of \( \Gamma \), thereby preserving the property that every edge corresponds with a proper hyperplane, or
  \item we require that edges between equal points have gain zero, thereby ensuring that each such edge corresponds with a degenerate hyperplane, or
  \item we make no additional requirements, instead adding caveats concerning equations of the form \( d^2(x, q) - d^2(x, q) = r \).
  \end{enumerate}
  As none of these adds anything of present interest, our assumption is harmless.
\end{rmk}

After coordinatizing our affine space, we may express the left hand side of the hyperplane equations in a Pythagorean arrangement as follows using the Euclidean norm.
\begin{align*}
  d^2(x, q_u) - d^2(x, q_v)
  & = |x - q_u|^2 - |x - q_v|^2 \\
  & = (x - q_u) \cdot (x - q_u) - (x - q_v) \cdot (x - q_v) \\
  & = (x \cdot x - 2 x \cdot q_u + q_u \cdot q_u) - (x \cdot x - 2 x \cdot q_v + q_v \cdot q_v) \\
  & = 2x \cdot (q_v - q_u) - (|q_v|^2 - |q_u|^2)
    .
\end{align*}
Thus the equations of the hyperplane \( h_e \) may be rewritten as follows:
\begin{align*}
  & d^2(x,q_u) - d^2(x, q_v) = g(e) \\
  \iff
  & 2x \cdot (q_v - q_u) - (|q_v|^2 - |q_u|^2) = g(e) \\
  \iff
  & 2x \cdot (q_v - q_u) = g(e) + |q_v|^2 - |q_u|^2
    .
\end{align*}

\begin{ex}
  \label{example:affinographic-case}
  Given a gain graph \( (\Gamma, g) \), the \emph{affinographic arrangement} of this data is a hyperplane arrangement in \( \R^{V(\Gamma)} \) with one hyperplane for each edge.
  For an edge \( e \) oriented from \( u \) to \( v \), the associated hyperplane has equation \( x_v - x_u = g(e) \).
  We recover affinographic arrangements as a special class of Pythagorean arrangements by choosing appropriate reference points.
  Treating \( \R^{V(\Gamma)} \) as an affine space, let \( q_v = \tfrac{1}{2}\mathbf{1}_v \) for each \( v \in V(\Gamma) \), where \( \mathbf{1}_\alpha \) denotes the indicator vector of \( \alpha \).
  Let \( \mc{Q} = \set{q_v}{v \in V(\Gamma)} \).
  In the Pythagorean arrangement \( \pythag{\mc{Q}}{\Gamma}{g} \), the hyperplane associated to an edge \( e \) oriented from \( u \) to \( v \) has equation
  \begin{align*}
    g(e)
    & = d^2(x, q_u) - d^2(x, q_v) \\
    & = 2x \cdot(q_v - q_u) - (|q_v|^2 - |q_u|^2) \\
    & = 2x \cdot(\tfrac{1}{2}\mathbf{1}_v - \tfrac{1}{2}\mathbf{1}_u) - (|\tfrac{1}{2}\mathbf{1}_v|^2 - |\tfrac{1}{2}\mathbf{1}_u|^2) \\
    & = x \cdot (\mathbf{1}_v - \mathbf{1}_u) - \tfrac{1}{2}(1 - 1) \\
    & = x_v - x_u
    ,
  \end{align*}
  which is precisely the equation of the affinographic hyperplane for the edge \( e \).
\end{ex}

A subset \( S \subseteq E(\Gamma) \) is \emph{central} for \( (\mc{Q}, \Gamma, g) \) when it corresponds to a central subarrangement of \( \pythag{\mc{Q}}{\Gamma}{g} \).
Central circles are of particular interest in this work.
\begin{prop}
  \label{result:central-circles}
  Let \( (\mc{Q}, \Gamma, g) \) be given and let \( C \) be a circle in \( \Gamma \).
  \begin{enumerate}
  \item \label[part]{result:central-circles:central-implies-balanced} %
    If \( C \) is central, then \( C \) is balanced.
  \item \label[part]{result:central-circles:balanced-lin-indep-implies-central} %
    If \( C \) is balanced and there is an edge \( e \) of \( C \) such that
    \[
      \set{q_v - q_u}{C \setminus \{e\} \text{ has an edge from } u \text{ to } v }
    \]
    is linearly independent, then \( C \) is central.
  \item \label[part]{result:central-circles:balanced-short-lin-indep-implies-central} %
    If \( C \) is a balanced circle of length at most \( d+1 \) and the points of \( \mc{Q} \) corresponding to vertices of \( C \) are in affine position, then \( C \) is central.
  \end{enumerate}
\end{prop}
We need the following lemma; this is an elementary exercise in affine and projective geometry, rephrased to invoke graph theory.
\begin{lem}
  \label{result:ideal-position-d+1}
  Let \( \mc{Q} \) be a collection of \( d+1 \) points in \( \affineSpace^d \).
  The following are equivalent.
  \begin{enumerate}
  \item %
    The point set \( \mc{Q} \) is in ideal position.
  \item %
    For every (equivalently, some) oriented spanning tree \( T \) in the complete graph on vertex set \( \mc{Q} \), the set \( B_T \) is linearly independent in the direction space \( D(\affineSpace) \).
  \item %
    \( \dim(D(\affineSpan{\mc{Q}})) = d \)
  \item %
    \( \dim(\affineSpan{\mc{Q}}) = d \)
  \item %
    The point set \( \mc{Q} \) is in affine position.
  \end{enumerate}
  Furthermore, a point set \( \mc{Q} \) of arbitrary size \( n \) is in ideal position precisely when \( n \geq d+1 \) and every \( (d+1) \)-subset thereof is in ideal position.
  \noProof
\end{lem}

\begin{proof}[Proof of \cref{result:central-circles}]
  Express \( C \) as \( (e_1, e_2, \dots, e_k) \) where \( e_i \) is an edge from \( v_{i-1} \) to \( v_i \) with \( v_0 = v_k \).

  To prove \cref{result:central-circles:central-implies-balanced}, assume \( C \) is central and let \( x \) be an element of the intersection.
  We compute
  \begin{align*}
    \sum_{i=1}^k g(e)
    & = \sum_{i=1}^k (d^2(x, q_{i-1}) - d^2(x, q_{i})) \\
    & = \sum_{i=1}^k d^2(x, q_{i-1}) - \sum_{i=1}^k d^2(x, q_{i}) \\
    & = d^2(x, q_0) + \left(\sum_{i=2}^k d^2(x, q_{i-1}) - \sum_{i=1}^{k-1} d^2(x, q_{i})\right) - d^2(x, q_k) \\
    & = d^2(x, q_0) - d^2(x, q_0) + \sum_{i=1}^{k-1} (d^2(x, q_i) - d^2(x, q_{i})) \\
    & = 0
      ,
  \end{align*}
  which yields that \( C \) is balanced.

  To prove \cref{result:central-circles:balanced-lin-indep-implies-central}, let \( B \) be any basis of \( D(\affineSpace) \) containing the vectors \( b_i = q_i - q_{i-1} \) for \( 1 \leq i \leq k \).
  As \( B \) is a basis, the linear system with equations
  \[
    \left\{
      \begin{array}{rclccr}
        2 x \cdot (q_i - q_{i-1}) &=& g(e_i) + |q_i|^2 - |q_{i-1}|^2 &&& 1 \leq i \leq k \\
        2 x \cdot b_i            &=& 0                          &&& k < i \leq d
      \end{array}
    \right.
  \]
  is uniquely solvable.
  For each \( 1 \leq i \leq k \), the solution \( x \) of this linear system satisfies
  \begin{align*}
    d^2(x, q_{i-1}) - d^2(x, q_i)
    & = |x - q_{i-1}|^2 - |x - q_i|^2 \\
    & = (x - q_{i-1}) \cdot (x - q_{i-1}) - (x - q_i) \cdot (x - q_i) \\
    & = (x \cdot x - 2 x \cdot q_{i-1} + q_{i-1} \cdot  q_{i-1}) - (x \cdot x - 2 x \cdot q_i + q_i \cdot  q_i) \\
    & = 2 x \cdot (q_i - q_{i-1}) + |q_{i-1}|^2 - |q_i|^2 \\
    & = (g(e_i) + |q_i|^2 - |q_{i-1}|^2) + |q_{i-1}|^2 - |q_i|^2 \\
    & = g(e_i)
      .
  \end{align*}
  Hence \( x \) is on \( h_{e_i} \) for all \( i \), and thus \( C \) is central.

  \cref{result:central-circles:balanced-short-lin-indep-implies-central} follows from \cref{result:central-circles:balanced-lin-indep-implies-central} and \cref{result:ideal-position-d+1}.
\end{proof}
By the \emph{linear system associated to a subgraph \( \Delta \subseteq \Gamma \)}, we mean the linear system
\begin{align}
  \label{equation:associated-linear-system}
  \left\{d^2(x, q_u) - d^2(x, q_v) = g(e)\phantom{\bigg|}\right.
  &&
     e \text{ is an edge of } \Delta \text{ from } u \text{ to } v
     .
\end{align}
The following simple, but useful, proposition follows from our proof of \cref{result:central-circles:central-implies-balanced}.
\begin{sco}
  \label{result:overdetermined-central-circles}
  The linear system \eqref{equation:associated-linear-system} of any central circle is over-determined.
\end{sco}

\begin{cor}
  If \( \Delta \) is a balanced subgraph of \( \Gamma \) with rank \( r \leq d \) and the points of \( \mc{Q} \) are in affine position, then the linear system associated to \( \Delta \) is consistent.
\end{cor}
\begin{proof}
  Choose a forest \( F \subseteq \Delta \) with \( r \).
  Every edge \( e \in \Delta \setminus F \) determines a balanced circle with at most \( r+1 \) edges.
  Thus the equation associated to the edge \( e \) is a linear combination of the equations of \( F \) by \cref{result:overdetermined-central-circles}.
  Hence \eqref{equation:associated-linear-system} is consistent.
\end{proof}

\subsection{%
  Modular ideals and linear classes%
  \label{section:modular-ideals}%
}

A \emph{linear class of circuits} of \( M \) is a subset \( S \subseteq \circuits{M} \) with the property that for every modular family of circuits \( T \subseteq S \) and every \( X \in \circuits{M} \) with \( X \subseteq \bigcup T \) we have \( X \in S \).
We begin by proving that linear classes carry the same information as modular ideals, i.e., there is a natural isomorphism between their lattices.

\begin{rmk}
  Some of what follows in this section is implicit in \cite{biased_graphs_2__the_three_matroids--zaslavsky}, though Zaslavsky was interested in the special case for \emph{biased graphs}, i.e., graphs with a distinguished linear class of circles.
  The matroids case is noted there without proof.
\end{rmk}

\begin{prop}
  \label{result:linear-classes-are-closed-sets}
  The set of linear classes of circuits in \( M \) is the family of closed sets of a closure operator on \( \circuits{M} \).
\end{prop}
Recall that a set system \( \mc{F} \subseteq \powerset{E} \) is the system of closed sets of some closure operator precisely when \( E \in \mc{F} \) and for all \( X, Y \in \mc{F} \) we have \( X \cap Y \in \mc{F} \).
\begin{proof}
  Note that \( \circuits{M} \) is trivially a linear class of circuits.
  If \( \linearClassSymbol_1 \) and \( \linearClassSymbol_2 \) are linear classes, then for all modular pairs \( X, Y \in \linearClassSymbol_1 \cap \linearClassSymbol_2 \), every circuit contained in \( X \cup Y \) belongs to both \( \linearClassSymbol_1 \) and \( \linearClassSymbol_2 \) by assumption.
  Thus \( \linearClassSymbol_1 \cap \linearClassSymbol_2 \) is a linear class by definition.
\end{proof}

\begin{figure}
  \centering
  \newcommand{\makeGraphNodes}[1]{
    \draw[scale=1.25, every node/.style={draw,circle, inner sep=1.5pt, fill=black}]
    (1,1) node (1) {}
    (0,2) node (2) {}
    (0,0) node (3) {}
    (2,0) node (4) {}
    (2,2) node (5) {}
    ;
    \tikzset{mine/.style={ultra thick, #1}}
  }

  \begin{subfigure}{.3\textwidth}
    \centering
    \begin{tikzpicture}
      \makeGraphNodes{Green}
      \draw[every node/.style={circle, inner sep=2pt, fill=white, draw=none}]
      (1) edge[mine] node {$a$} (2) edge node {$b$} (3) edge[mine] node {$c$} (4)
      (2) edge[mine] node {$d$} (3) edge node {$e$} (5)
      (3) edge[mine] node {$f$} (4)
      (4) edge node {$g$} (5)
      ;
    \end{tikzpicture}
    \subcaption{\( \textcolor{Green}{X} \subseteq E(\Gamma) \)}
  \end{subfigure}
  \begin{subfigure}{.3\textwidth}
    \centering
    \begin{tikzpicture}
      \makeGraphNodes{blue}
      \draw[every node/.style={circle, inner sep=2pt, fill=white, draw=none}]
      (1) edge[mine] node {$a$} (2) edge[mine] node {$b$} (3) edge node {$c$} (4)
      (2) edge node {$d$} (3) edge[mine] node {$e$} (5)
      (3) edge[mine] node {$f$} (4)
      (4) edge[mine] node {$g$} (5)
      ;
    \end{tikzpicture}
    \subcaption{\( \textcolor{blue}{Y} \subseteq E(\Gamma) \)}
  \end{subfigure}
  \begin{subfigure}{.3\textwidth}
    \centering
    \begin{tikzpicture}
      \makeGraphNodes{red}
      \draw[every node/.style={circle, inner sep=2pt, fill=white, draw=none}]
      (1) edge node {$a$} (2) edge[mine] node {$b$} (3) edge[mine] node {$c$} (4)
      (2) edge[mine] node {$d$} (3) edge node {$e$} (5)
      (3) edge[mine] node {$f$} (4)
      (4) edge[mine] node {$g$} (5)
      ;
    \end{tikzpicture}
    \subcaption{\( \textcolor{red}{A} \subseteq E(\Gamma) \)}
  \end{subfigure}
  \caption{%
    \label{figure:modular-nonclosure-graphs}%
    Modular ideal closures are not idempotent.
  }
\end{figure}
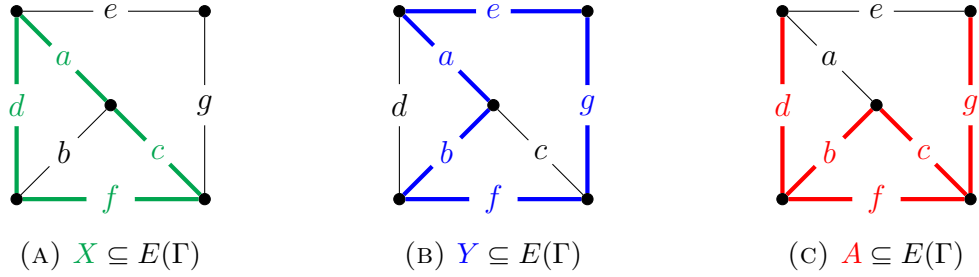

We call the lattice of linear classes of \( M \) the \emph{modular circuit lattice}.
\begin{ex}
  One might be tempted to define a closure analogous to the usual matroid closure on \( E \) using a linear class \( S \subseteq \mc{C} \) of circuits as follows for all \( A \subseteq E \).
  \[
    \sigma A
    = A \cup \set{e \in E}{e \in X \subseteq A \cup \{e\} \text{ for some }X \in S}
    .
  \]
  While extensive and isotone, this operator need not be idempotent.
  In \cref{figure:modular-nonclosure-graphs}, the circles \( X = adfc \) and \( Y = aegfb \) together form a linear class \( S = \{X, Y\} \), but for \( A = bcdfg \) we have \( \sigma \sigma A = \sigma A \cup Y \neq \sigma A \).
\end{ex}

An ideal \( \mc{I} \subseteq \powerset{E} \) satisfies \emph{modular extension by circuits} when for all \( X \in \mc{I} \cap \circuits{M} \), all \( S \in \mc{I} \), and all \( e \in E \) with \( e \in X \subseteq S \cup \{e\} \) we have \( S \cup \{e\} \in \mc{I} \).
\begin{thm}
  \label{result:linear-classes-yield-modular-ideals}
  Let \( M \) be a matroid with order ideal \( \mc{I} \subseteq \powerset{E} \).
  The following are equivalent.
  \begin{enumerate}
  \item %
    \( \mc{I} \) is a modular ideal of \( M \).
  \item\label{ideal-via-circuits} %
    \begin{enumerate}
    \item %
      \( \mc{I} \cap \circuits{M} \) is a linear class of circuits in \( M \),
    \item %
      \( \independentSets{M} \subseteq \mc{I} \), and
    \item %
      \( \mc{I} \) satisfies modular extension by circuits.
    \end{enumerate}
  \end{enumerate}
\end{thm}
\begin{proof}
  Suppose \( X, Y \) form a modular pair of circuits in \( \mc{I} \cap \circuits{M} \).
  Then \( X, Y \) form a modular pair in \( \mc{I} \), so \( X \cup Y \in \mc{I} \).
  As \( \mc{I} \) is an ideal, every circuit \( Z \in \circuits{M} \) with \( Z \subseteq X \cup Y \) belongs to \( \mc{I} \).
  Hence \( Z \in \mc{I} \cap \circuits{M} \), and \( \mc{I} \cap \circuits{M} \) is a linear class.
  We noted previously that the independent sets all belong to \( \mc{I} \).
  Finally, assume \( X \in \mc{I} \cap \circuits{M} \), \( S \in \mc{I} \), and \( e \in E \) satisfy \( e \in X \subseteq S \cup \{e\} \).
  If \( e \in S \), the conclusion holds trivially.
  Otherwise, \( \rk(S \cap X) = \cardinality{X}-1 = \rk(X) \) and \( \rk(S \cup X) = \rk(S) \), so \( \rk(S \cap X) + \rk(S \cup X) = \rk(S) + \rk(X) \); thus these form a modular pair and \( S \cup \{e\} = S \cup X \in \mc{I} \) as desired.

  Suppose \( \mc{I} \cap \circuits{M} \) is a linear class of circuits in \( M \), \( \independentSets{M} \subseteq \mc{I} \), and that we may apply modular extension by circuits.
  Let \( S, T \in \mc{I} \) be an arbitrary modular pair.
  Let \( B \) be any basis for \( S \cup T \) obtained by extending a basis for \( S \cap T \) first to a basis for \( S \) and then to a basis for \( S \cup T \).
  By construction \( \cardinality{B} = \rk(S \cup T) \), \( \cardinality{B \cap S \cap T} = \rk(S \cap T) \), and \( \cardinality{B \cap S} = \rk(S) \); by modularity
  \[
    \cardinality{B}
    = \rk(S \cup T)
    = \rk(S) + \rk(T) - \rk(S \cap T)
    = \cardinality{B \cap S} - \cardinality{B \cap S \cap T} + \rk(T)
    .
  \]
  Thus \( \rk(T) = \cardinality{B \cap T} \) by inclusion-exclusion, so \( B \cap T \) is a basis of \( T \).

  Now we prove \( S \cup T \in \mc{I} \).
  Order \( (S \cup T) \setminus B = \{x_1, \dots, x_m\} \) and define \( U_k = B \cup \{x_1, \dots, x_k\} \).
  As \( B \) is independent, we see \( U_0 = B \in \mc{I} \) by assumption.
  Assume that \( U_{k-1} \in \mc{I} \) for some \( k \in [m] \).
  If \( x_k \in S \), then \( C_B(x_k) \subseteq (B \cap S) \cup \{x_k\} \subseteq S \) because \( B \cap S \) is a basis of \( S \); similarly, if \( x_k \in T \), then \( C_B(x_k) \subseteq T \).
  In either case, \( C_B(x_k) \subseteq I \in \mc{I} \) for either \( I = S \) or \( I = T \), so \( C_B(x_k) \in \mc{I} \) as \( \mc{I} \) is an order ideal.
  Finally,
  \[
    x_k
    \in C_B(x_k)
    \subseteq B \cup \{x_k\}
    \subseteq U_{k-1} \cup \{x_k\}
    ,
  \]
  so we may apply modular extension by circuits to obtain \( U_k = U_{k-1} \cup \{x_k\} \in \mc{I} \).
  Hence \( S \cup T = U_m \in \mc{I} \) by induction, and \( \mc{I} \) satisfies the modular extension property.
  As non-degeneracy is implied by \( \independentSets{M} \subseteq \mc{I} \), we have \( \mc{I} \) is a modular ideal.
\end{proof}

Given a basis \( B \) of \( S \subseteq E \) we define \( F_B(S) = \set{C_B(s)}{s \in S \setminus B} \).
\begin{cor}
  The elements of a modular ideal \( \mc{I} \) are precisely
  \[
    \independentSets{M}
    \cup
    \set{S \subseteq E}
    {\text{there is a basis \( B \) of \( S \) such that } F_B(S) \subseteq \mc{I}}
    .
  \]
\end{cor}
\begin{proof}
  If \( S \in \mc{I} \), then every element of \( F_B(S) \) is a member of \( \mc{I} \) as this is an order ideal.
  If there is a basis \( B \) of \( S \) for which \( F_B(S) \subseteq \mc{I} \), then \( B \in \mc{I} \) and thus repeated application of extension by circuits yields \( S \in \mc{I} \).
\end{proof}

\begin{cor}
  The lattice of modular ideals of \( M \) is isomorphic to the lattice of linear classes of circuits in \( M \).
\end{cor}
\begin{proof}
  Note \( f(\mc{I}) = \mc{I} \cap \circuits{M} \) preserves the lattice meet, which is intersection in both cases, and trivially sends \( \top \) to \( \top \).
  Moreover, if \( f(\mc{I}_1) = f(\mc{I}_2) \), then \( \mc{I}_1 = \mc{I}_2 \) by the previous corollary.
  Finally, given \( \linearClassSymbol \subseteq \circuits{M} \) a linear class of circuits, let
  \[
    \mc{I}_{\linearClassSymbol}
    = \set{S \subseteq E}
    {\exists B \subseteq S \text{ with } \rk(S) = \cardinality{B} \text{ and } F_B(S) \subseteq \linearClassSymbol}
    .
  \]
  Thus \( \mc{I}_{\linearClassSymbol} \cap \circuits{M} = \linearClassSymbol \) by the previous corollary.
  Hence \( f \) is a meet semilattice isomorphism preserving \( \top \), and it is thus a lattice isomorphism.
\end{proof}

\subsection{Combinatorics of Pythagorean arrangements%
  \label{section:pythagorean-combinatorics}%
}

Fix a point set \( \mc{Q} \) in affine \( d \)-space and a graph \( \Gamma \) whose vertices are in bijection with \( \mc{Q} \); the point of \( \mc{Q} \) associated to a vertex \( v \) of \( \Gamma \) is denoted \( q_v \).
We think of a gain function \( g \) on \( \Gamma \) as an element of the \emph{(real) edge space} of \( \Gamma \), i.e., the elements \( x \in \R^{\vec{E}(\Gamma)} \) satisfying \( x(\bar{e}) = -x(e) \) for all oriented edges \( e \) of \( \Gamma \), where \( \bar{e} \) denotes \( e \) with the opposite orientation.
This space is linearly isomorphic with \( \R^{E(\Gamma)} \); indeed, choosing a preferred orientation \( \omega \) of the edges of \( \Gamma \) yields a linear map \( L_\omega \colon \R^{E(\Gamma)} \to \R^{\vec{E}(\Gamma)} \) defined by \( L_\omega (\mathbf{1}_e) = \mathbf{1}_{e^+} - \mathbf{1}_{e^-} \), where \( \mathbf{1}_\alpha \) is the indicator function of \( \alpha \).
This map \( L_\omega \) sends a basis of \( \R^{E(\Gamma)} \) to a basis of edge space.
Going forward, we work in \( \R^{E(\Gamma)} \) via one such orientation, which we mostly ignore.
This choice does not significantly affect results, though it does make statements cleaner.

A triple \( (\mc{Q}, \Gamma, g) \) is \emph{gain-generic} when there is an \( \epsilon > 0 \) such that for all gain functions \( g' \), if the Euclidean norm in edge space satisfies \( |g' - g| < \epsilon \), then the Pythagorean arrangements \( \pythag{\mc{Q}}{\Gamma}{g'} \) and \( \pythag{\mc{Q}}{\Gamma}{g} \) have the same combinatorial type, i.e., all such perturbations preserve the labeled intersection semilattice.

We now describe failures of gain genericity.
First we introduce a dummy variable \( \gamma_e \) for each edge \( e \) of \( \Gamma \), a placeholder for the gain on that edge.
Let \( S \) be a \( (d+1) \)-subset of \( E(\Gamma) \).
This yields a square matrix \( \systemMatrix{\Gamma}{\mc{Q}}{S} \) with \( d+1 \) rows and entries in the polynomial ring \( \R[\gamma_e : e \in E(\Gamma)] \) as follows.
Each edge \( e \in S \) yields a row:\ the first \( d \) entries are the entries of \( q_v - q_u \) where \( e \) is oriented from \( u \) to \( v \), and the final entry is \( -\tfrac{1}{2}(\gamma_e + |q_v|^2 - |q_u|^2) \).

Suppose the variables are all specified via some gain function, i.e., set \( \gamma_e = g(e) \) for edges \( e \), and let \( S \subseteq E(\Gamma) \) have rank \( d \) in \( M_\Gamma(Q) \).
Define \( \systemMatrix{\Gamma}{Q}{S} = [ A \ -b] \).
By elementary linear algebra, the over-determined linear system \( Ax = b \) has a solution precisely when \( \det(\systemMatrix{\Gamma}{Q}{S}) = 0 \).
Furthermore, the corresponding linear system is uniquely solvable, i.e., the \( d+1 \) hyperplanes associated to the elements of \( S \) intersect in a unique common point because the codimension of their intersection is \( d \).
Moreover this is non-generic:\ we may vary any one of the involved gains \( g(e) \) by an arbitrarily small \( \epsilon > 0 \) to obtain a new linear system.
As the original was uniquely solvable and over-determined, the new system has no solution.

In terms of the corresponding Pythagorean arrangements, the original gains corresponded with a central subarrangement of \( d+1 \) hyperplanes, having exactly one point in the intersection.
After varying one gain, the corresponding subarrangement is non-central.
This process effectively nudged one of the hyperplanes off of the unique central point by an arbitrarily small change of gain function.
Hence \( \det(\systemMatrix{\Gamma}{Q}{S}) = 0 \) yields a failure of gain genericity.

For each circuit of \( M_\Gamma(\mc{Q}) \) expressed as a fundamental circuit with respect to some basis, the above construction witnesses a failure of genericity, i.e., if the gain function satisfies the equation \( \det(\systemMatrix{\Gamma}{Q}{S}) = 0 \), then the arrangement is non-generic, and the circuit is gain-generically non-central.
Our main result in this section is that the converse also holds, i.e., that every failure of gain genericity has a corresponding central circuit \( C \) of \( M_\Gamma(\mc{Q}) \).

Note that the equation \( \det(\systemMatrix{\Gamma}{Q}{S}) = 0 \) resulting from our argument above is linear in the gain variables; this equation thus yields a linear hyperplane in edge space.

\begin{figure}
  \newcommand{\gCD}{6}
  \newcommand{\gAB}{-\gCD}
  \newcommand{\gAC}{0}
  \newcommand{\gAD}{2}
  \newcommand{\gBC}{\gAD}

  \begin{subfigure}[b]{.35\textwidth}
    \centering
    \begin{tikzpicture}
      \draw[every node/.style={circle, fill, inner sep=1.5pt}]
      (0,0) node[label={225:$1$}] (a) {}
      (2,0) node[label={315:$2$}] (b) {}
      (2,2) node[label={45:$3$}] (c) {}
      (0,2) node[label={135:$4$}] (d) {}
      ;
      \draw[
      every node/.style={inner sep=1.5pt},
      >=Stealth,
      arr/.style={decoration={markings, mark={at position .5 with {\arrow{>}}}}, postaction={decorate}}
      ]
      (a) edge[arr] node[label={-90:$\gAB$}] {} (b)
      (a) edge[arr] node[label={135:$\gAC$}] {} (c)
      (a) edge[arr] node[label={180:$\gAD$}] {} (d)
      (b) edge[arr] node[label={0:$\gBC$}] {} (c)
      (c) edge[arr] node[label={90:$\gCD$}] {} (d)
      ;
    \end{tikzpicture}
    \subcaption{Gain graph \( (\Gamma, g) \).}
  \end{subfigure}
  \begin{subfigure}[b]{.6\textwidth}
    \centering
    \begin{tikzpicture}
      \draw
      (0,0) node[circle, fill, inner sep=1pt, label={180:$q_1$}] (a) {}
      (4,0) node[circle, fill, inner sep=1pt, label={0:$q_2$}] (b) {}
      (3,2) node[circle, fill, inner sep=1pt, label={90:$q_3$}] (c) {}
      (1,2) node[circle, fill, inner sep=1pt, label={90:$q_4$}] (d) {}
      ;
      \clip(-1,-1) rectangle (5,3);
      \pgfmathsetmacro{\xAB}{2+\gAB/8}
      \pgfmathsetmacro{\xCD}{2-\gCD/4}
      \draw [thick] (\xAB,-1) -- (\xAB,3);
      \draw [thick,domain=-1:5] plot(\x,{(13+\gAC-6*\x)/4});
      \draw [thick,domain=-1:5] plot(\x,{(5+\gAD-2*\x)/4});
      \draw [thick,domain=-1:5] plot(\x,{(\gBC-3+2*\x)/4});
      \draw [thick] (\xCD,-1) -- (\xCD,3);
    \end{tikzpicture}
    \subcaption{Arrangement \( \pythag{\mc{Q}}{\Gamma}{g} \)}
  \end{subfigure}
  \caption{%
    \label{figure:gain-generic-but-balanced}%
    Gain-generic Pythagorean arrangement with a balanced circle.
  }
\end{figure}
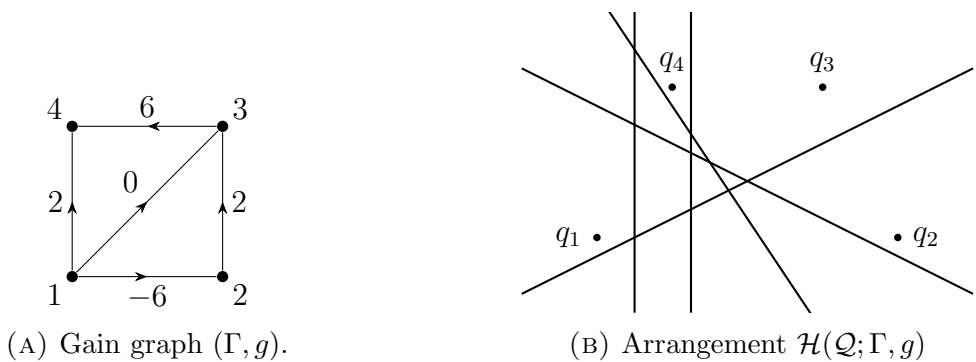

\begin{ex}
  \label{example:forbidden-gains}
  Consider the graph and point configuration illustrated in \cref{figure:graph-and-configuration}.
  \begin{figure}
    \begin{subfigure}[b]{.35\textwidth}
      \centering
      \begin{tikzpicture}
        \draw[every node/.style={circle, fill, inner sep=1.5pt}]
        (0,0) node[label={225:$1$}] (a) {}
        (2,0) node[label={315:$2$}] (b) {}
        (2,2) node[label={45:$3$}] (c) {}
        (0,2) node[label={135:$4$}] (d) {}
        ;
        \draw[
        every node/.style={inner sep=1.5pt},
        >=Stealth,
        arr/.style={decoration={markings, mark={at position .5 with {\arrow{>}}}}, postaction={decorate}}
        ]
        (a) edge[arr] node[label={-90:$a$}] {} (b)
        (a) edge[arr] node[label={135:$b$}] {} (c)
        (a) edge[arr] node[label={180:$c$}] {} (d)
        (b) edge[arr] node[label={  0:$s$}] {} (c)
        (c) edge[arr] node[label={ 90:$t$}] {} (d)
        ;
      \end{tikzpicture}
      \subcaption{Graph \( \Gamma \).}
    \end{subfigure}
    \begin{subfigure}[b]{.6\textwidth}
      \centering
      \begin{tikzpicture}
        \draw
        (0,0) node[circle, fill, inner sep=1pt, label={-90:$q_1=(0,0)$}] (a) {}
        (4,0) node[circle, fill, inner sep=1pt, label={-90:$q_2=(4,0)$}] (b) {}
        (3,2) node[circle, fill, inner sep=1pt, label={0:$q_3=(3,2)$}] (c) {}
        (1,2) node[circle, fill, inner sep=1pt, label={180:$q_4=(1,2)$}] (d) {}
        ;
      \end{tikzpicture}
      \subcaption{Configuration \( \mc{Q} \).}
    \end{subfigure}
    \caption{%
      \label{figure:graph-and-configuration}%
      Graph and point configuration for \cref{example:forbidden-gains}.
    }
  \end{figure}
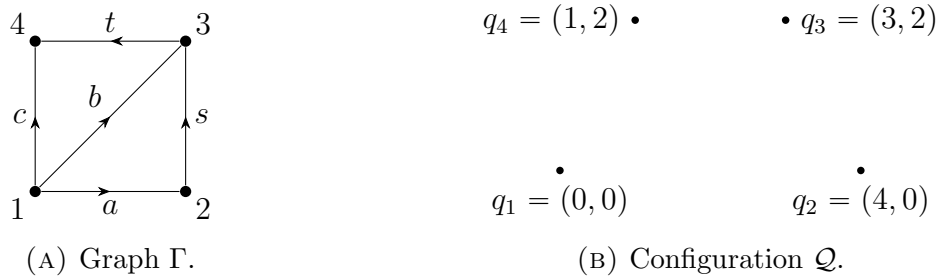
  Let \( M = M_\infty(\mc{Q})|E(\Gamma) \).
  Each triple of edges in \( \Gamma \) has a corresponding \emph{hyperplane of nongenericity} in edge space, \( \det(\systemMatrix{\Gamma}{Q}{S}) = 0 \), obtained from the preceding discussion.
  These are presented in the table below with some simplifications by canceling common constant factors.
  \[
    \begin{array}{crccr}
      \text{Edges}~~&\text{Relation} &~~~~~~~~~& \text{Edges}~~&\text{Relation} \\
      \hline
      abc & \gamma_a - 2\gamma_b + 2\gamma_c = 0 & & ast & \gamma_a + 2\gamma_t = 0         \\
      abs & \gamma_a - \gamma_b + \gamma_s = 0   & & bcs & \gamma_b - 2\gamma_c + \gamma_s = 0   \\
      abt & \gamma_a + 2\gamma_t = 0             & & bct & \gamma_b - \gamma_c + \gamma_t = 0    \\
      acs & \gamma_a - 2\gamma_c + 2\gamma_s = 0 & & bst & \gamma_b - \gamma_s + 2\gamma_t = 0   \\
      act & \gamma_a + 2\gamma_t = 0             & & cst & \gamma_c - \gamma_s + \gamma_t = 0
    \end{array}
  \]

  The equations corresponding to the sets \( abs \) and \( bct \) are the equations of the ``short circles'' in the graph \( \Gamma \), i.e., circles with at most \( d+1 = 3\) edges.
  Each such equation is precisely the equation we would expect to hold if the circle were balanced, precisely because each short circle has its corresponding reference points in ideal position for this configuration.

  The equation \( \gamma_a + 2\gamma_t = 0 \) results from edge sets \( abt \), \( act \), and \( ast \).
  The support of this equation is \( at \), a circuit of \( M_\Gamma(\mc{Q}) \) from a projective dependence in our point set---indeed, there is a parallelism \( \affineSpan{q_1, q_2} \parallel \affineSpan{q_3, q_4} \).
  Removing either \( a \) or \( t \) from each triad yields a basis of \( M \) with \( at \) a fundamental circuit, e.g., \( at = C_{ac}(t) = C_{ct}(a) \).

  Finally, the circle with edges \( astc \) is irrelevant in this construction.
  Indeed, each \( k \)-gon in \( \Gamma \) with \( k > d+1 = 3 \) contains a circuit of \( M_\Gamma(\mc{Q}) \) with at most \( d+1 \) edges.
  Our construction detects failures of genericity on these smaller sets.
  In this case, if \( astc \) were central, we would see that the circuits \( at \), \( acs \), and \( cst \) were also central.
  \cref{result:pythagorean-arrangement-combinatorics} proves \( S \subseteq E(\Gamma) \) is central if and only if every circuit \( C \subseteq S \) of \( M_\Gamma(\mc{Q}) \) is central.
  Thus these ``long circles'' with at least \( d+2 = 4 \) edges do not correspond to a hyperplane of nongenericity.
  \cref{figure:gain-generic-but-balanced} illustrates a gain-generic arrangement for \( \mc{Q} \), with balanced long circle \( astc \).
\end{ex}

\begin{rmk}
  Our discussion of long circles is significant because of how starkly this contrasts the affinographic case of \cref{example:affinographic-case}.
  In such cases, there are no long circles due to the restriction on dimension.
  Indeed, the simplex configuration \( \mc{Q} \) is in ideal position; thus there are no circuits of \( M_\Gamma(\mc{Q}) \) resulting from projective dependence.
  Hence, the equations resulting from these circuits ought all be balancing equations for the corresponding circles of \( \Gamma \).
  In fact, \( M_\Gamma(\mc{Q}) \cong M[\Gamma] \).

  In light of this example, our main result will recover the classical result on affinographic arrangements that the semimatroid of the affinographic arrangement of \( (\Gamma, g) \) is isomorphic to the semimatroid of balanced subgraphs of \( (\Gamma, g) \).
\end{rmk}

Let \( M_\Gamma(\mc{Q}) = M_\infty(\mc{Q})|E(\Gamma) \), i.e., \( M_\Gamma(\mc{Q}) \) is the matroid at infinity of \( \mc{Q} \) restricted to elements corresponding to edges of \( \Gamma \).
For each circuit \( X \) of \( M_\Gamma(\mc{Q}) \), choose an \( x \in X \) and a basis \( B \) of \( M_\Gamma(\mc{Q}) \) for which \( X = C_B(x) \).
From this data, define
\begin{equation}
  \label{equation:forbidden-hyperplane-definition}
  F_X
  =
  \set{g \in \R^{E(\Gamma)}}
  {\det(\systemMatrix{\Gamma}{\mc{Q}}{B \cup \{x\}}) = 0}
  .
\end{equation}
We call this the \emph{hyperplane of nongenericity} for \( X \).
\begin{lem}
  \label{result:gain-generic-forbidden-hyperplanes-are-well-defined}
  The hyperplane \( F_X \) is well-defined, i.e., \( F_X \) is independent of the choices of \( x \in X \) and basis \( B \) made.
  Moreover, the equation of \( F_X \) is supported on the variable symbols \( \gamma_e \) with \( e \in X \).
\end{lem}
Below we conflate a vector \( q_e \) in the ground set of \( M_\Gamma(\mc{Q}) \) with the corresponding edge \( e \) of \( \Gamma \) to simplify our language.
Recall that the linear system associated to a subgraph \( \Delta \subseteq \Gamma \) is
\begin{align*}
  \left\{d^2(x, q_u) - d^2(x, q_v) = g(e)\phantom{\bigg|}\right.
  &&
     e \text{ is an edge of } \Delta \text{ from } u \text{ to } v
     \tag{\ref{equation:associated-linear-system}}
     ,
\end{align*}
which can be rewritten as
\begin{align*}
  \left\{2(q_v - q_u) \cdot x = g(e) + |q_v|^2 - |q_u|^2\phantom{\bigg|}\right.
  &&
     e \text{ is an edge of } \Delta \text{ from } u \text{ to } v
     \tag{\ref{equation:associated-linear-system}'}
     .
\end{align*}
\begin{proof}[Proof of \cref{result:gain-generic-forbidden-hyperplanes-are-well-defined}]
  We first prove that for every \( x \in X \) and every basis \( B = \{b_1, b_2, \dots, b_d\} \) of \( M_\Gamma(\mc{Q}) \) with \( X \setminus \{x\} \subseteq B \), the polynomial \( \det (\systemMatrix{\Gamma}{\mc{Q}}{B \cup x}) \) is supported on the symbols \( \gamma_e \) with \( e \in X \).
  Let \( b_i \) be oriented from \( u_i \) to \( v_i \) for each \( i \in \{0\} \cup [d] \) where \( b_0 = x \), and recall that \( b_i = q_{v_i} - q_{u_i} \).
  Consider the matrix
  \[
    \systemMatrix{\Gamma}{\mc{Q}}{B \cup \{x\}}
    =
    \begin{bmatrix}
      q_{v_0} - q_{u_0} & -\tfrac{1}{2}(\gamma_{b_0} + |q_{v_0}|^2 - |q_{u_0}|^2) \\
      q_{v_1} - q_{u_1} & -\tfrac{1}{2}(\gamma_{b_1} + |q_{v_1}|^2 - |q_{u_1}|^2) \\
      q_{v_2} - q_{u_2} & -\tfrac{1}{2}(\gamma_{b_2} + |q_{v_2}|^2 - |q_{u_2}|^2) \\
      \vdots & \vdots \\
      q_{v_d} - q_{u_d} & -\tfrac{1}{2}(\gamma_{b_d} + |q_{v_d}|^2 - |q_{u_d}|^2)
    \end{bmatrix}
    .
  \]
  Let \( D_i \) denote the submatrix obtained from \( \systemMatrix{\Gamma}{\mc{Q}}{B \cup \{x\}} \) by removing the \( i^{th} \) row (indexed from \( 0 \)) and final column.
  That is,
  \[
    D_i
    =
    \begin{bmatrix}
      q_{v_0} - q_{u_0} \\
      \vdots \\
      q_{v_{i-1}} - q_{u_{i-1}} \\
      q_{v_{i+1}} - q_{u_{i+1}} \\
      \vdots \\
      q_{v_d} - q_{u_d}
    \end{bmatrix}
    .
  \]
  The cofactor expansion of \( \det(\systemMatrix{\Gamma}{\mc{Q}}{B \cup \{x\}}) \) along the final column is thus
  \begin{align*}
    \det(\systemMatrix{\Gamma}{\mc{Q}}{B \cup \{x\}})
    & = (-1)^d \sum_{i=0}^d (-1)^i -\tfrac{1}{2}(\gamma_{b_i} + |q_{v_i}|^2 - |q_{u_i}|^2)\det D_i \\
    & = (-1)^{d+1} \tfrac{1}{2} \sum_{i=0}^d (-1)^i \det D_i \cdot (\gamma_{b_i} + |q_{v_i}|^2 - |q_{u_i}|^2)
      .
  \end{align*}
  By uniqueness of fundamental circuits, \( b_i \in X \) if and only if \( B_{\hat{i}} = (\{x\} \cup B) \setminus \{b_i\} \) is independent.
  On the other hand, the rows of \( D_i \) are precisely the vectors in \( B_{\hat{i}} \).
  Hence \( \det D_i \neq 0 \) if and only if \( b_i \in X \) by basic properties of the determinant.
  Finally, \( \tfrac{1}{2}(-1)^{d+i+1}\det D_i \) is the coefficient of the variable \( \gamma_{b_i} \) for all \( i \), so \( \gamma_{b_i} \) is in the support of the polynomial \( \det(\systemMatrix{\Gamma}{\mc{Q}}{B \cup \{x\}}) \) precisely when \( b_i \in X \).

  Next we show that for all \( x, x' \in X \) and all bases \( B \) and \( B' \) of \( M_\Gamma(\mc{Q}) \) for which \( C_B(x) = X = C_{B'}(x') \), the equations
  \begin{align*}
    \det(\systemMatrix{\Gamma}{\mc{Q}}{B \cup \{x\}}) &= 0
    &&\text{and}
    & \det(\systemMatrix{\Gamma}{\mc{Q}}{B' \cup \{x'\}}) &= 0
  \end{align*}
  are equivalent, i.e., have the same solution set.
  As the elements of \( M_\Gamma(\mc{Q}) \) are vectors in \( \R^d \) and \( B, B' \) are bases of \( \R^d \) by construction of \( M_\Gamma(\mc{Q}) \), there is a linear isomorphism \( L \) of \( \R^d \) for which \( B' = L(B) \).
  Now \( x = \sum_{b \in B} \alpha_b b \), so
  \[
    L(x)
    = \sum_{b \in B} \alpha_b L(b)
    = \sum_{b' \in B'} \alpha_b b'
    = x'
    .
  \]
  Thus when we compare the polynomials
  \begin{align*}
    \det(\systemMatrix{\Gamma}{\mc{Q}}{B \cup \{x\}})
    && \text{and}
    && \det(\systemMatrix{\Gamma}{\mc{Q}}{B' \cup \{x'\}})
       ,
  \end{align*}
  we are comparing the matrices below.
  \begin{align*}
    \systemMatrix{\Gamma}{\mc{Q}}{B \cup \{x\}}
    & =
      \begin{bmatrix}
        q_{v_0} - q_{u_0} & -\tfrac{1}{2}(\gamma_{b_0} + |q_{v_0}|^2 - |q_{u_0}|^2) \\
        q_{v_1} - q_{u_1} & -\tfrac{1}{2}(\gamma_{b_1} + |q_{v_1}|^2 - |q_{u_1}|^2) \\
        q_{v_2} - q_{u_2} & -\tfrac{1}{2}(\gamma_{b_2} + |q_{v_2}|^2 - |q_{u_2}|^2) \\
        \vdots & \vdots \\
        q_{v_d} - q_{u_d} & -\tfrac{1}{2}(\gamma_{b_d} + |q_{v_d}|^2 - |q_{u_d}|^2)
      \end{bmatrix}
    \\
    \systemMatrix{\Gamma}{\mc{Q}}{B' \cup \{x'\}}
    & =
      \begin{bmatrix}
        L(q_{v_0} - q_{u_0}) & -\tfrac{1}{2}(\gamma_{b_0} + |q_{v_0}|^2 - |q_{u_0}|^2) \\
        L(q_{v_1} - q_{u_1}) & -\tfrac{1}{2}(\gamma_{b_1} + |q_{v_1}|^2 - |q_{u_1}|^2) \\
        L(q_{v_2} - q_{u_2}) & -\tfrac{1}{2}(\gamma_{b_2} + |q_{v_2}|^2 - |q_{u_2}|^2) \\
        \vdots & \vdots \\
        L(q_{v_d} - q_{u_d}) & -\tfrac{1}{2}(\gamma_{b_d} + |q_{v_d}|^2 - |q_{u_d}|^2)
      \end{bmatrix}
  \end{align*}
  Preserving the notation from before, the submatrix
  \[
    D_i'
    =
    \begin{bmatrix}
      q_{v'_0} - q_{u'_0} \\
      \vdots \\
      q_{v'_{i-1}} - q_{u'_{i-1}} \\
      q_{v'_{i+1}} - q_{u'_{i+1}} \\
      \vdots \\
      q_{v'_d} - q_{u'_d}
    \end{bmatrix}
    =
    \begin{bmatrix}
      L(q_{v_0} - q_{u_0}) \\
      \vdots \\
      L(q_{v_{i-1}} - q_{u_{i-1}}) \\
      L(q_{v_{i+1}} - q_{u_{i+1}}) \\
      \vdots \\
      L(q_{v_d} - q_{u_d})
    \end{bmatrix}
    =
    \begin{bmatrix}
      q_{v_0} - q_{u_0} \\
      \vdots \\
      q_{v_{i-1}} - q_{u_{i-1}} \\
      q_{v_{i+1}} - q_{u_{i+1}} \\
      \vdots \\
      q_{v_d} - q_{u_d}
    \end{bmatrix}
    \cdot
    [L]^t
    =
    D_i \cdot [L]^t
    ,
  \]
  where \( [L]^t \) denotes the transpose of the matrix representing the linear map \( L \) in the standard basis on \( \R^{E(\Gamma)} \).
  Thus, by basic properties of the determinant,
  \[
    \det(D'_i)
    = \det(D_i \cdot [L]^t)
    = \det(D_i) \det(L)
    .
  \]
  Hence the cofactor expansion of \( \det(\systemMatrix{\Gamma}{\mc{Q}}{B' \cup \{x'\}}) \) along the final column is
  \begin{align*}
    \det(\systemMatrix{\Gamma}{\mc{Q}}{B' \cup \{x'\}})
    & = (-1)^{d+1} \tfrac{1}{2} \sum_{i=0}^d (-1)^i \det D'_i \cdot (\gamma_{b'_i} + |q_{v'_i}|^2 - |q_{u'_i}|^2) \\
    & = (-1)^{d+1} \tfrac{1}{2} \det(L) \sum_{i=0}^d (-1)^i \det D_i \cdot (\gamma_{b_i} + |q_{v_i}|^2 - |q_{u_i}|^2) \\
    & = \det(L) \det(\systemMatrix{\Gamma}{\mc{Q}}{B \cup \{x\}})
      .
  \end{align*}
  Finally, note that \( \det(L) \neq 0 \) as \( L \) is an isomorphism.
  Hence the equations corresponding to \( B \cup \{x\} \) and \( B' \cup \{x'\} \) are equivalent.
\end{proof}

The most important tool to prove our main theorem is the arrangement consisting of all hyperplanes \( F_X \) described above, i.e.,
\[
  \mc{F}(\mc{Q}, \Gamma)
  = \set{F_X}{X \in \circuits{M_\Gamma(\mc{Q})}}
  .
\]
First, we establish a strong connection between gain-genericity and this arrangement.
\begin{prop}
  \label{result:gain-generic-forbidden-equations}
  The triple \( (\mc{Q}, \Gamma, g) \) is gain-generic if and only if \( g \notin F_X \) for all \( X \in \circuits{M_\Gamma(\mc{Q})} \).
\end{prop}
\begin{proof}
  We first show that each \( g \in F_X \) has non-generic \( (\mc{Q}, \Gamma, g) \).
  Setting \( \gamma_e = |q_v|^2 - |q_u|^2 \) for all \( e \) oriented from \( u \) to \( v \) results in a linear system of \( d+1 \) homogeneous equations; thus the system admits a solution, and the equation of \( F_X \) is satisfied.
  By multilinearity of the determinant, we may subtract this solution from the final column without affecting the determinant.
  The resulting column is in the span of the others precisely when the corresponding linear system (of \( d+1 \) equations in \( d \) unknowns) is solvable, i.e., when the determinant of our matrix is \( 0 \); but this matrix is precisely \( \systemMatrix{\Gamma}{\mc{Q}}{B \cup X} \) where \( X \) is a fundamental circuit with respect to \( B \).
  Thus the \( d+1 \) hyperplanes of \( B \cup X \) form a central subarrangement of \( \pythag{\mc{Q}}{\Gamma}{g} \) and \( (\mc{Q}, \Gamma, g) \) is non-generic.

  Conversely, assume \( (\mc{Q}, \Gamma, g) \) is not gain-generic.
  Thus there is a central subarrangement of \( \pythag{\mc{Q}}{\Gamma}{g} \) with \( d+1 \) hyperplanes.
  The corresponding hyperplane directions form a dependent set \( D \) of \( M_\Gamma(\mc{Q}) \), which in turn contains a circuit \( X \).
  Choose any basis \( B \) for which \( X \) is a fundamental circuit; then assigning \( \gamma_e = g(e) \) for all edges \( e \) yields \( \det(\systemMatrix{\Gamma}{\mc{Q}}{B \cup X}) = 0 \), so \( g \in F_X \).
  Hence every non-generic \( g \) belongs to some hyperplane of \( \mc{F}(\mc{Q}, \Gamma) \).
\end{proof}

\begin{cor}
  If \( (Q, \Gamma, g) \) is gain-generic, then the circles of \( (\Gamma, g) \) with at most \( d+1 \) edges are all unbalanced.
\end{cor}
\begin{proof}
  Every such circle contains a circuit of \( M_\Gamma(\mc{Q}) \).
\end{proof}

A subset of the ground set of a matroid is \emph{unicyclic} when it contains a unique circuit.
In light of the above corollary, a unicyclic set of \( M_\Gamma(\mc{Q}) \) contains at most one circle of \( \Gamma \), and cannot have more than \( d+1 \) elements.
The maximal unicyclic subsets always have the form \( B \cup \{x\} \) for some basis \( B \) and \( x \in E \setminus B \).
Note that these unicyclic sets are precisely the sets used to define the hyperplanes \( F_X \).
\begin{sco}
  For every intersection flat of \( \pythag{\mc{Q}}{\Gamma}{g} \), there is either an independent or unicyclic set of \( M_\Gamma(\mc{Q}) \) yielding the same intersection.
\end{sco}
Combining this with the well-definition portion of the argument yields the following.
\begin{sco}
  \label{result:central-iff-on-forbidden-hyperplane}
  A gain function \( g \) belongs to \( F_X \) if and only if \( X \) is a central circuit of \( \pythag{\mc{Q}}{\Gamma}{g} \).
\end{sco}

Given \( g \in \R^{E(\Gamma)} \) we define the flat of \( g \) by
\[
  A_g
  = \bigcap \set{F_X}{X \in \circuits{M_\Gamma(\mc{Q})} \text{ and } g \in F_X}
  .
\]

\begin{thm}
  \label{result:pythagorean-arrangement-combinatorics}
  Let \( (\mc{Q}, \Gamma) \) be given.
  \begin{enumerate}
  \item %
    For each flat \( A \) of \( \mc{F}(\mc{Q}, \Gamma) \) the set
    \[
      \circuits{A}
      = \set{X \in \circuits{M_\Gamma(\mc{Q})}}{F_X \supseteq A}
    \]
    is a linear class of circuits in \( M_\Gamma(\mc{Q}) \).
  \item %
    For each \( g \) in edge space, the set of central circuits of \( \pythag{\mc{Q}}{\Gamma}{g} \) is \( \circuits{A_g} \).
  \item %
    The lattice of flats of \( \mc{F}(\mc{Q}, \Gamma) \) is isomorphic to the lattice of semimatroids realizable as \( \pythag{\mc{Q}}{\Gamma}{g} \) for some gain function \( g \).
  \end{enumerate}
\end{thm}
Informally, the intersection pattern of \( \pythag{\mc{Q}}{\Gamma}{g} \) is determined entirely by the corresponding flat of \( g \) in the arrangement \( \mc{F}(\mc{Q}, \Gamma) \).
Moreover, the hyperplanes determining each flat contain the information necessary to reconstruct the central subarrangements of all \( \pythag{\mc{Q}}{\Gamma}{g} \) with gain function \( g \) belonging to this flat.
The proof is simple, albeit somewhat technical, using following fundamental lemma from linear algebra.
\begin{lem}[Three-Term Pl\"{u}cker Relations]
  \label{lemma:plucker-relations}
  Let \( E = \{e_1, e_2, \dots, e_{n-2}\} \subseteq \R^n \), and for each \( x, y \in \R^n \) let \( [x, y; E] \) denote the determinant of the matrix whose rows are the vectors \( x, y, e_1, e_2, \dots, e_n \) in this order.
  For all \( a, b, c, d \) we have
  \[
    [a, b; E][c, d; E] - [a, c; E][b, d; E] + [a, d; E][b, c; E]
    = 0
    .
    \noProof
  \]
\end{lem}
\begin{proof}[Proof of \cref{result:pythagorean-arrangement-combinatorics}]
  Suppose \( X, Y \in \mc{I}(A) \) are a modular pair and \( Z \subseteq X \cup Y \); if \( Z = X \) or \( Z = Y \), then \( Z \in \mc{C}(A) \) by well-definition of \( F_Z \) in \cref{result:gain-generic-forbidden-equations}.
  Thus we may assume \( X \), \( Y \), and \( Z \) are distinct.
  Let \( B \) be a basis of \( M \) with \( X = C_B(x) \) and \( Y = C_B(y) \).
  Now \( X \cap Y \neq \emptyset \) lest eliminating \( x \) between \( X \) and \( Z \) yield a circuit
  \[
    Z'
    \subseteq (X \cup Z) \setminus \{x\}
    \subseteq (X \cup Y) \setminus \{x\}
    \subsetneq X \cup Y
    ,
  \]
  contradicting modularity of \( X \) and \( Y \).
  Moreover, if \( X \cap Y \subseteq Z \), then eliminating \( z \in X \cap Y \) between the pairs \( X, Z \) and \( Y, Z \) yields circuits \( X' \) and \( Y' \) with
  \begin{align*}
    X' \cup Y'
    & \subseteq ((X \cup Z) \setminus \{z\}) \cup ((Y \cup Z) \setminus \{z\}) \\
    & = (X \cup Y \cup Z) \setminus \{z\} \\
    & = (X \cup Y) \setminus \{z\} \\
    & \subsetneq X \cup Y
      ,
  \end{align*}
  which again contradicts modularity of the pair \( X, Y \).
  Hence \( X \cap Y \setminus Z \neq \emptyset \).

  Let \( z \in X \cap Y \setminus Z \) and note that the circuit elimination of \( z \) between \( X \) and \( Y \) is necessarily \( Z \); otherwise we have a circuit \( Z' \subseteq (X \cup Y) \setminus \{z\} \), and this yields \( Z \cup Z' \subseteq (X \cup Y) \setminus \{z\} \subsetneq X \cup Y \), again contradicting the assumption that \( X, Y \) is a modular pair.
  Consider the \( (d+2)\times(d+1) \) matrix
  \[
    W
    =
    \begin{bmatrix}
      q_x         & G_x     \\
      q_y         & G_y     \\
      q_z         & G_z     \\
      q_{b_1}     & G_{b_1} \\
      q_{b_2}     & G_{b_2} \\
      \vdots      & \vdots  \\
      q_{b_{d-1}} & G_{b_{d-1}}
    \end{bmatrix}
  \]
  where \( q_e \) denotes the (row) vector \( q_v - q_u \) for edge \( e \) from \( u \) to \( v \), \( G_e = -\tfrac{1}{2}(\gamma_e + q_v - q_u) \), and \( B' = \{b_1, b_2, \dots, b_{d-1}\} = B \setminus \{z\} \).
  Omitting the row corresponding to an edge \( e \in \{x, y, z, b_1, \dots, b_{d-1}\} \) yields the matrix \( \systemMatrix{\mc{Q}}{\Gamma}{(B' \cup \{x,y,z\}) \setminus \{e\}} \), i.e., one of
  \begin{align*}
    W_Y
    & =
      \begin{bmatrix}
        q_y         & G_y     \\
        q_z         & G_z     \\
        q_{b_1}     & G_{b_1} \\
        q_{b_2}     & G_{b_2} \\
        \vdots      & \vdots  \\
        q_{b_{d-1}} & G_{b_{d-1}}
      \end{bmatrix}
                      ,
                    &
                      W_X
                    & =
                      \begin{bmatrix}
                        q_x         & G_x     \\
                        q_z         & G_z     \\
                        q_{b_1}     & G_{b_1} \\
                        q_{b_2}     & G_{b_2} \\
                        \vdots      & \vdots  \\
                        q_{b_{d-1}} & G_{b_{d-1}}
                      \end{bmatrix}
                                      ,
                    &
                      W_{Z}
                    & =
                      \begin{bmatrix}
                        q_x         & G_x     \\
                        q_y         & G_y     \\
                        q_{b_1}     & G_{b_1} \\
                        q_{b_2}     & G_{b_2} \\
                        \vdots      & \vdots  \\
                        q_{b_{d-1}} & G_{b_{d-1}}
                      \end{bmatrix}
                                      ,
  \end{align*}
  or
  \[
    W_{b_k}
    =
    \begin{bmatrix}
      q_x         & G_x     \\
      q_y         & G_y     \\
      q_z         & G_z     \\
      q_{b_1}     & G_{b_1} \\
      q_{b_2}     & G_{b_2} \\
      \vdots      & \vdots  \\
      q_{b_{k-1}} & G_{b_2} \\
      q_{b_{k+1}} & G_{b_2} \\
      \vdots      & \vdots  \\
      q_{b_{d-1}} & G_{b_{d-1}}
    \end{bmatrix}
  \]
  Recall the defining equations of \( F_\alpha \) from \eqref{equation:forbidden-hyperplane-definition}.
  Thus \( F_\alpha = \set{g}{\det(W_\alpha) = 0} \) for each \( \alpha = X, Y, Z \).
  With the notation of \cref{lemma:plucker-relations}, for each \( k \in [d-1] \) the Pl\"{u}cker relation in our matrix corresponding to \( y, b_k, x, z \) and \( B'_{\hat{k}} = B' \setminus \{b_k\} \) is
  \[
    0
    =
    [y, b_k; B'_{\hat{k}}][x, z; B'_{\hat{k}}]
    - [y, x; B'_{\hat{k}}][b_k, z; B'_{\hat{k}}]
    + [y, z; B'_{\hat{k}}][b_k, x; B'_{\hat{k}}]
    .
  \]
  Note \( B = B' \cup \{z\} \), so solving for the middle term yields
  \begin{align*}
    (-1)^k[z;B'][x, y; B'_{\hat{k}}]
    & =
      (-1)^{k+1}[z; B'][y, x; B'_{\hat{k}}]
    \\
    & =
      [b_k, z; B'_{\hat{k}}][y, x; B'_{\hat{k}}]
    \\
    & =
      [y, b_k; B'_{\hat{k}}][x, z; B'_{\hat{k}}]
      + [b_k, x; B'_{\hat{k}}][y, z; B'_{\hat{k}}]
    \\
    & =
      (-1)^k[y; B'][x, z; B'_{\hat{k}}]
      + (-1)^{k+1}[x; B'][y, z; B'_{\hat{k}}]
    \\
    & =
      (-1)^{k+1}
      \left(
      [x; B'][y, z; B'_{\hat{k}}]
      - [y; B'][x, z; B'_{\hat{k}}]
      \right)
      .
  \end{align*}
  Dividing by \( (-1)^{k+1} \) yields
  \[
    [z; B'][x, y; B'_{\hat{k}}]
    = [y; B'][x, z; B'_{\hat{k}}] - [x; B'][y, z; B'_{\hat{k}}]
    ,
  \]
  and the coefficients of the \( [-,-;B'_{\hat{k}}] \)'s are independent of \( k \).
  Moreover, the left side is precisely the coefficient of \( G_{b_k} \) as in the definition for \( F_Z \) in \eqref{equation:forbidden-hyperplane-definition}, and the right side is a linear combination of the coefficients of \( G_{b_k} \) in defining equations for \( F_X \) and \( F_Y \).
  Denoting the \( (\beta, d+1) \)-cofactor of \( W_\alpha \) by \( W_\alpha(\beta) \) for \( \alpha = X, Y, Z \), we have shown
  \[
    [z;B']W_Z(b_k) = [y;B']W_X(b_k) - [x;B']W_Y(b_k)
  \]
  for all \( k \in [d-1] \).
  The same coefficients work with \( G_z \):
  \begin{align*}
    [z;B']W_Z(z)
    & = [z;B'] 0
      = 0 \\
    & = [y;B'] 0 - [x;B'] 0 \\
    & = [y;B'][z,z;B'] - [x;B'][z,z;B'] \\
    & = [y;B']W_X(z) - [x;B']W_Y(z)
      .
  \end{align*}
  They similarly work for \( G_x \) and \( G_y \):
  \begin{align*}
    [z;B']W_Z(x)
    & = [y;B'][z;B'] - [x;B'] 0
      = [y;B']W_X(x) - [x;B']W_Y(x)
      ,
    \\
    [z;B']W_Z(y)
    & = [z;B'](-[x;B'])
      = [y;B'] 0 - [x;B'][x;B'] \\
    & = [y;B']W_X(y) - [x;B']W_Y(y)
      .
  \end{align*}
  Moreover, the coefficient of \( G_e \) is \( 0 \) for all \( e \in E \setminus (B \cup \{x,y\}) \) in all equations for \( F_X \), \( F_Y \), and \( F_Z \) by \cref{result:gain-generic-forbidden-hyperplanes-are-well-defined}.
  Thus the following computation proves that the defining polynomial for \( F_Z \) is a linear combination of the polynomials defining \( F_X \) and \( F_Y \):
  \begin{align*}
    [z;B']\det(W_Z)
    & = \sum_{\beta \in B \cup \{x, y\}} [z; B']W_Z(\beta)G_\beta \\
    & = \sum_{\beta \in B \cup \{x, y\}} ([y;B']W_X(\beta) - [x;B']W_Y(\beta))G_\beta \\
    & = \sum_{\beta \in B \cup \{x, y\}} [y;B']W_X(\beta)G_\beta
      - \sum_{\beta \in B \cup \{x, y\}} [x;B']W_Y(\beta)G_\beta \\
    & = [y;B']\det(W_X) - [x;B']\det(W_Y)
  \end{align*}
  Thus \( F_X \cap F_Y \subseteq F_Z \); indeed, if \( \det(W_X) = 0 = \det(W_Y) \) for a particular gain function \( g \), then \( \det(W_Z) = 0 \).
  Hence \( Z \in \circuits{A} \), and \( \circuits{A} \) is a linear class of circuits.

  By the work above and \cref{result:central-iff-on-forbidden-hyperplane}, we have for all gain functions \( g \):
  \begin{align*}
    X \in \circuits{A_g}
    & \iff A_g \subseteq F_X \\
    & \iff X \text{ is a circuit with } \set{h_e}{e \in X} \text{ central in } \pythag{\mc{Q}}{\Gamma}{g}.
  \end{align*}
  Hence we have a complete characterization of the central circuits; he remainder of the work for the second part is done by \cref{result:linear-classes-yield-modular-ideals} and \cref{lemma:semimatroid-iff-modular-ideal}.

  The third part follows from the first two.
  The desired map takes the intersection semilattice of \( \pythag{\mc{Q}}{\Gamma}{g} \) to the flat \( A_g \).
  By the second part, the intersection semilattice of \( \pythag{\mc{Q}}{\Gamma}{g} \) is determined by \( \circuits{A_g} \), so this map is well-defined.
  A poset is realizable as an intersection semilattice of a Pythagorean arrangement if and only if there is a gain function realizing it, so this map is trivially surjective.
  Moreover, two flats of \( \mc{F}(\mc{Q}, \Gamma) \) are equal if and only if they are determined by the same maximal set of circuits of \( M_\Gamma(\mc{Q}) \), so this map is injective.
  Finally, for flats \( A, B \) of \( \mc{F}(\Gamma, \mc{Q}) \) we have \( A \subseteq B \) if and only if \( \circuits{A} \supseteq \circuits{B} \); recalling that the flats of a hyperplane arrangement are ordered by reverse inclusion, we conclude that this map is a poset isomorphism.
\end{proof}

We now amend \cref{result:pythagorean-arrangement-combinatorics} to deal with prescribed lists of balanced circles, i.e., with prescribed bias.
Fix a graph \( \Gamma \), a linear class of circles \( \linearClassSymbol \subseteq \circuits{\Gamma} \), and a point configuration \( \mc{Q} \) labeled by \( V(\Gamma) \).
Each circle \( X \in \circuits{\Gamma} \) yields a hyperplane in edge space, namely
\[
  B_X
  = \set{g \in \R^E}
  {\sum_{e \in X} g(e) = 0}
\]
where the edges are oriented coherently along the circle.
Note that this hyperplane may or may not belong to \( \mc{F}(\mc{Q}, \Gamma) \), depending on the geometry of \( \mc{Q} \).
Indeed, if any proper subset of \( X \) satisfies a nontrivial projective dependence relation, then \( X \) is not a circuit of \( M_\Gamma(\mc{Q}) \), and \( B_X \) is not an \( F_Y \).
Thus the ``long circles'' of \( \Gamma \), i.e., those with length \( l > d+1 \), determine a \( B_X \) but not an \( F_X \).

The set of gain functions satisfying exactly the balance prescribed by \( \linearClassSymbol \) is
\[
  \mathfrak{G}_{\linearClassSymbol}
  = \bigcap_{X \in \linearClassSymbol} B_X
  \setminus
  \bigcup_{Y \in \mc{C}(\Gamma) \setminus \linearClassSymbol} B_Y
  .
\]
This set may be empty if \( \linearClassSymbol \) is not realizable.
In the case that \( \linearClassSymbol \) is realizable, these observations allow us to conclude the following.
\begin{cor}
  The intersection patterns of Pythagorean arrangements on \( (\mc{Q}, \Gamma) \) with balanced circles exactly \( \linearClassSymbol \) are in bijection with the set
  \[
    \set{\mathfrak{G}_{\linearClassSymbol} \cap A}{A \in \flats{\mc{F}(\mc{Q}, \Gamma)} \text{ and } \mathfrak{G}_{\linearClassSymbol} \cap A \neq \emptyset}
    .
  \]
  The set of all such combinatorial types is a meet-semilattice ordered by weak map.
  \noProof
\end{cor}

Now we explore the structure imposed by edge space a bit more.
If \( B \) is a basis of \( M_\Gamma(\mc{Q}) \), then we may treat the gains on edges of \( B \) as free variables.
Then every fundamental circuit with respect to this basis has its gain determined if the circuit is central (again, the associated linear system is over-determined).
As not all circles are central, the general situation is quite complicated, and there may be relations that go undetected by this basis---we did need the whole hyperplane arrangement \( \mc{F}(\mc{Q}, \Gamma) \), after all.
On the other hand, if \( \pythag{\mc{Q}}{\Gamma}{g} \) is central, then the flat \( A_g \) is completely determined by the chosen basis \( B \).
Indeed, the fact that \( \bigcap_{e \in B} h_e \) is a point and the arrangement is central together determine all the remaining gains.
Fixing a basis \( B \) of \( M_\Gamma(\mc{Q}) \) determines a bijection between edge space and the ambient affine space.
Thus we obtain the following.
\begin{prop}
  \label{result:centres-correspond-to-gains-in-top-element}
  The points of affine space correspond with the elements of gain space in \( \top \in \mc{F}(\mc{Q}, \Gamma) \).
  Each basis of \( M_\Gamma(\mc{Q}) \) determines one such correspondence.
  \noProof
\end{prop}

\section{Examples and Applications}
\label{section:applications}

This section contains a variety of examples and applications of our results above.
It also includes some refinements and corollaries of these results for those applications.
All edges in examples are oriented from the smaller vertex to the greater, unless otherwise indicated.

\subsection{Small Example with Additional Bias}
{
  \newcommand{\scaleParam}{.7}
  \newcommand{\biasColor}{Red}
  \newcommand{\normalColor}{Blue}
  \newcommand{\gainGraph}[6]{%
    \begin{tikzpicture}%
      \draw[every node/.style={circle, fill, inner sep=1.5pt}]%
      (0,0) node (a) {}%
      (2,0) node (b) {}%
      (2,2) node (c) {}%
      (0,2) node (d) {}%
      ;%
      \draw[%
      every node/.style={inner sep=1.5pt},%
      >=Stealth,%
      arr/.style={decoration={markings, mark={at position .5 with {\arrow{>}}}}, postaction={decorate}}%
      ]%
      (a) edge[arr] node[below] {$#1$} (b)%
      (a) edge[arr] node[sloped,above] {$#2$} (c)%
      (a) edge[arr] node[sloped,above] {$#3$} (d)%
      (b) edge[arr] node[right] {$#4$} (c)%
      (c) edge[arr] node[above] {$#5$} (d)%
      ;%
      \ifthenelse{\equal{#6}{}}{}{\node[draw=Green, inner sep=1pt, scale=1.5] at (1,-.75) {$#6$};}
    \end{tikzpicture}%
  }%
  \newcommand{\cAB}{Green}
  \newcommand{\cAC}{blue}
  \newcommand{\cAD}{Red}
  \newcommand{\cBC}{Orange}
  \newcommand{\cCD}{Purple}
  \newcommand{\myArrangement}[7][scale=.4]{%
    \begin{tikzpicture}[
      >=latex,
      hyp/.style={domain=-1:5},
      hypEdge/.style={ultra thin, dashed},
      arr/.style={dashed},
      transform shape,
      #1
      ]
      \node[draw=Green, inner sep=1pt, scale=1.5] at (2,-.75) {$#2$};
      \clip (-.25,-.25) rectangle (4.25,2.25);
      \pgfmathsetmacro{\xA}{2+#3/8}
      \pgfmathsetmacro{\gB}{#4}
      \pgfmathsetmacro{\gC}{#5}
      \pgfmathsetmacro{\gS}{#6}
      \pgfmathsetmacro{\xT}{2-#7/4}
      \path[every node/.style={circle, fill, inner sep=1.5pt}]
      (0,0) node (a) {}
      (4,0) node (b) {}
      (3,2) node (c) {}
      (1,2) node (d) {}
      ;
      \draw
      (a) edge[hypEdge,\cAB] (b)
      (a) edge[hypEdge,\cAC] (c)
      (a) edge[hypEdge,\cAD] (d)
      (b) edge[hypEdge,\cBC] (c)
      (c) edge[hypEdge,\cCD] (d)
      ;
      \draw[hyp, \cAB] (\xA,-1) -- (\xA,3);
      \draw[hyp, \cAC] plot(\x,{(13+\gB-6*\x)/4});
      \draw[hyp, \cAD] plot(\x,{(5+\gC-2*\x)/4});
      \draw[hyp, \cBC] plot(\x,{(\gS-3+2*\x)/4});
      \draw[hyp, \cCD] (\xT,-1) -- (\xT,3);
    \end{tikzpicture}
  }
  Consider the point configuration and graph of \cref{example:forbidden-gains}.
  The matroid \( M_\Gamma(\mc{Q}) \) has the modular circuits lattice depicted in \cref{figure:mod-circs-gain-lattice}, which is  isomorphic to the lattice \( L \) of flats of \( \mc{F}(\mc{Q}, \Gamma) \).
  \begin{figure}
    \centering
    \begin{tikzpicture}[scale=.95, every node/.style={inner sep=2pt, scale=\scaleParam, Blue}]
      \draw (0,10) node (1) {\gainGraph{-2t}   {s-2t}{s-t}{s}{t}{STUVWXYZ}};
      \draw[shift={(0,7)}, xscale=3.5]
      (-2,0) node (stuw)    {\gainGraph{2(c-s)}{2c-s}{c}  {s}{t}{STUW}}
      (-1,0) node (tvy)     {\gainGraph{-2t}   {s-2t}{c}  {s}{t}{TVY}}
      ( 0,0) node (uvx)     {\gainGraph{-2t}   {b}   {s-t}{s}{t}{UVX}}
      ( 1,0) node (svz)     {\gainGraph{-2t}   {c-t} {c}  {s}{t}{SVZ}}
      ( 2,0) node (wxyz)    {\gainGraph{a}     {s-2t}{s-t}{s}{t}{WXYZ}}
      ;
      \draw[shift={(0,3)}, xscale=2]
      (-3.5,0) node (s)     {\gainGraph{2(b-c)}{b}   {c}  {s}{t}{S=abc}}
      (-2.5,0) node (t)     {\gainGraph{b-s}   {b}   {c}  {s}{t}{T=abs}}
      (-1.5,0) node (u)     {\gainGraph{2(c-s)}{b}   {c}  {s}{t}{U=acs}}
      ( -.5,0) node (v)     {\gainGraph{-2t}   {b}   {c}  {s}{t}{V=at}}
      (  .5,0) node (w)     {\gainGraph{a}     {2c-s}{c}  {s}{t}{W=bcs}}
      ( 1.5,0) node (x)     {\gainGraph{a}     {b}   {s-t}{s}{t}{X=cst}}
      ( 2.5,0) node (y)     {\gainGraph{a}     {s-2t}{c}  {s}{t}{Y=bst}}
      ( 3.5,0) node (z)     {\gainGraph{a}     {c-t} {c}  {s}{t}{Z=bct}}
      ;
      \draw
      (0,0) node (0)        {\gainGraph{a}     {b}   {c}  {s}{t}{\emptyset}}
      ;
      \draw
      (1) edge (svz.north) edge (tvy.north) edge (uvx.north) edge (stuw.north) edge (wxyz.north)
      (svz.south) edge (s.north) edge (v.north) edge (z.north)
      (tvy.south) edge (t.north) edge (v.north) edge (y.north)
      (uvx.south) edge (u.north) edge (v.north) edge (x.north)
      (stuw.south) edge (s.north) edge (t.north) edge (u.north) edge (w.north)
      (wxyz.south) edge (w.north) edge (x.north) edge (y.north) edge (z.north)
      (0) edge (s.south) edge (t.south) edge (u.south) edge (v.south) edge (w.south) edge (x.south) edge (y.south) edge (z.south)
      ;
    \end{tikzpicture}
    \caption[Lattice of flats of \( \mc{F}(\mc{Q}, \Gamma) \).]{%
      \label{figure:mod-circs-gain-lattice}%
      Lattice of flats of \( \mc{F}(\mc{Q}, \Gamma) \).
      Lattice elements are labeled by a gain graph representing the flat and the modular family of circuits describing it.
      Circuit abbreviations label the atoms:\ uppercase letters are circuits, lowercase letters are edges.
    }
  \end{figure}
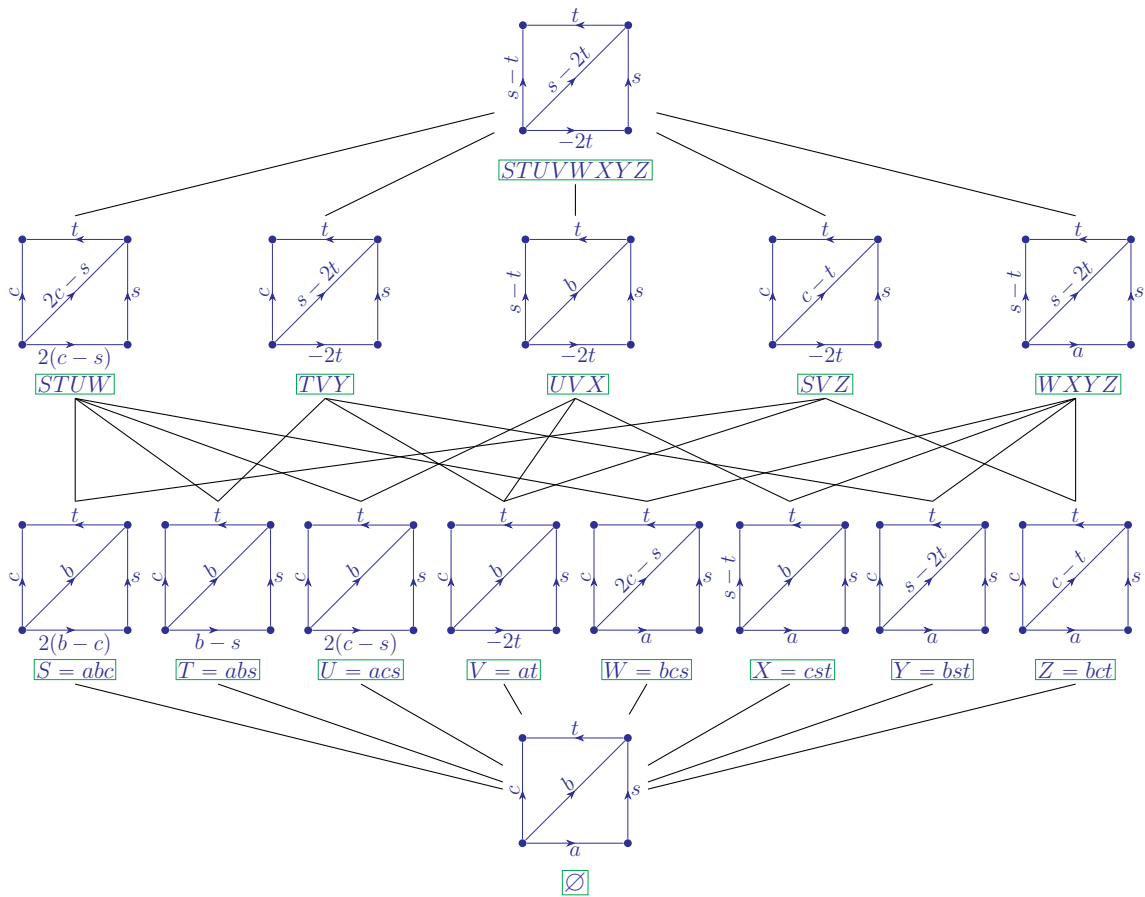
  As with every such arrangement, the bottom element \( \bot \) is the ambient space.

  The top element of \( L \) is a flat of dimension \( 2 \) in edge space (recall that \( L \) is ordered by reverse inclusion).
  Note that \( \top \) corresponds to the collection of all central arrangements \( \pythag{\mc{Q}}{\Gamma}{g} \).
  The modular ideal of central sets in this case is \( \powerset{E} \), so every subset is central.
  Thus if \( (\mc{Q}, \Gamma) \) is a full-dimensional pair in \( \affineSpace^d \) (i.e., if the directions spanned by the vectors \( q_e = q_v - q_u \) for edges \( e \) of \( \Gamma \) has dimension \( d \)), then so is the top element.
  As above, there is a bijection between \( \top \) and possible central points for each basis of \( M_\Gamma(\mc{Q}) \).
  \cref{figure:lattice-of-arrangements} has an example arrangement associated to each flat of \( \mc{F}(\mc{Q}, \Gamma) \).
  \begin{figure}
    \centering
    \begin{tikzpicture}[scale=.95] %
      \path
      (0,10)   node (1)     {\myArrangement{STUVWXYZ}{-2*0}{2-2*0}{2-0}{2}{0}}
      ;
      \draw[shift={(0,7)}, xscale=3.5]
      (-2,0)   node (stuw) {\myArrangement{STUW}{2*(2-0)}{2*2-0}{2}{0}{0}}
      (-1,0)   node (tvy)  {\myArrangement{TVY}{-2*0}{1-2*0}{0}{1}{0}}
      ( 0,0)   node (uvx)  {\myArrangement{UVX}{-2*(-1)}{1}{1-(-1)}{1}{-1}}
      ( 1,0)   node (svz)  {\myArrangement{SVZ}{-2*1}{2-1}{2}{1}{1}}
      ( 2,0)   node (wxyz) {\myArrangement{WXYZ}{0}{1-2*1}{1-1}{1}{1}}
      ;
      \draw[shift={(0,3)}, xscale=2]
      (-3.5,0) node (s)    {\myArrangement{S=abc}{2*(2-5)}{2}{5}{0}{1}}
      (-2.5,0) node (t)    {\myArrangement{T=abs}{0-4}{0}{1}{4}{0}}
      (-1.5,0) node (u)    {\myArrangement{U=acs}{2*(2-6)}{0}{2}{6}{1}}
      ( -.5,0) node (v)    {\myArrangement{V=at}{-2*1}{3}{3}{0}{1}}
      (  .5,0) node (w)    {\myArrangement{W=bcs}{1}{2*2-0}{2}{0}{1}}
      ( 1.5,0) node (x)    {\myArrangement{X=cst}{-3}{1}{5-0}{5}{0}}
      ( 2.5,0) node (y)    {\myArrangement{Y=bst}{0}{3-2*1}{0}{3}{1}}
      ( 3.5,0) node (z)    {\myArrangement{Z=bct}{-6}{1-0}{1}{3}{0}}
      ;
      \path
      (0,0)    node (0)     {\myArrangement{\emptyset}{1}{2}{0}{5}{3}}
      ;
      \draw
      (1) edge (svz.north) edge (tvy.north) edge (uvx.north) edge (stuw.north) edge (wxyz.north)
      (svz.south) edge (s.north) edge (v.north) edge (z.north)
      (tvy.south) edge (t.north) edge (v.north) edge (y.north)
      (uvx.south) edge (u.north) edge (v.north) edge (x.north)
      (stuw.south) edge (s.north) edge (t.north) edge (u.north) edge (w.north)
      (wxyz.south) edge (w.north) edge (x.north) edge (y.north) edge (z.north)
      (0) edge (s.south) edge (t.south) edge (u.south) edge (v.south) edge (w.south) edge (x.south) edge (y.south) edge (z.south)
      ;
    \end{tikzpicture}
    \caption{%
      \label{figure:lattice-of-arrangements}%
      Lattice of example arrangements corresponding with \cref{figure:mod-circs-gain-lattice}.
    }
  \end{figure}
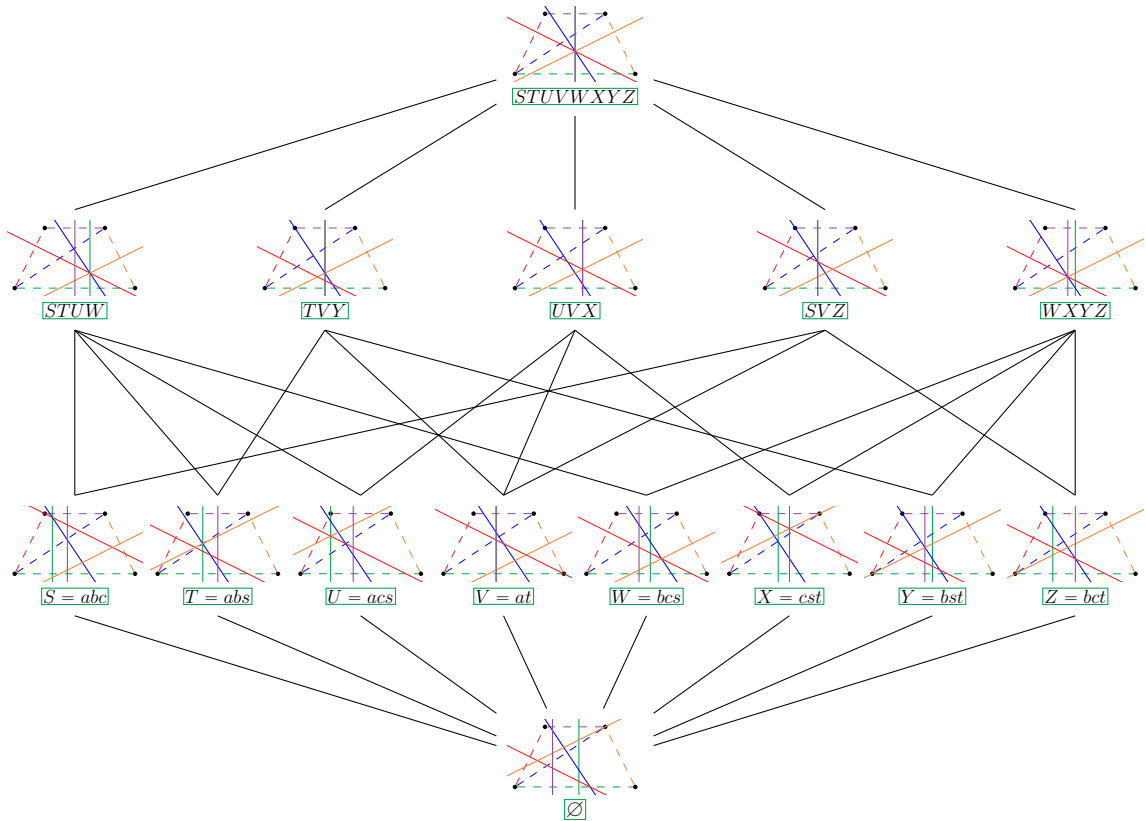

  Next we add bias; \( \Gamma \) is a theta graph, so our balanced circle classes must have either one or three circles.
  This gives four nonempty possibilities:\ \( \linearClassSymbol_1 = \{abd\} \), \( \linearClassSymbol_2 = \{bce\} \), \( \linearClassSymbol_3 = \{adec\} \), and \( \linearClassSymbol_4 = \{abd, bce, adec\} \).
  Enforcing these biases and preserving the labeling from \cref{figure:mod-circs-gain-lattice} yields restricted lattices of flats in \cref{figure:bias-1,figure:bias-2,figure:bias-3,figure:bias-4}.
  The flats which have more balanced circles than prescribed by \( \beta \) are highlighted in \MakeLowercase{\biasColor}.
  \begin{figure}
    \centering
    \begin{tikzpicture}[every node/.style={inner sep=2pt, scale=\scaleParam, \normalColor}]
      \draw (0,6) node[\biasColor] (1) {\gainGraph{-2t}   {s-2t}{s-t}{s}{t}{STUVWXYZ}};
      \draw[shift={(0,3)}, xscale=3]
      (-1.5,0) node (stuw)             {\gainGraph{2(c-s)}{2c-s}{c}  {s}{t}{STUW}}
      ( -.5,0) node (tvy)              {\gainGraph{-2t}   {s-2t}{c}  {s}{t}{TVY}}
      (  .5,0) node (tx)               {\gainGraph{b-s}   {b}   {s-t}{s}{t}{TX}}
      ( 1.5,0) node[\biasColor] (tz)   {\gainGraph{c-s-t} {c-t} {c}  {s}{t}{TZ}}
      ;
      \draw
      (0,0) node (0)                   {\gainGraph{b-s}   {b}   {c}  {s}{t}{T}}
      ;
      \draw
      (1) edge (stuw.north) edge (tvy.north) edge (tx.north) edge (tz.north)
      (0) edge (stuw.south) edge (tvy.south) edge (tx.south) edge (tz.south)
      ;
    \end{tikzpicture}
    \caption[Modular circuits with balance \( \linearClassSymbol_1 \).]{%
      \label{figure:bias-1}%
      \( \linearClassSymbol_1: a=b-s \).
    }
  \end{figure}

  \begin{figure}
    \centering
    \begin{tikzpicture}[every node/.style={inner sep=2pt, scale=\scaleParam, \normalColor}]
      \draw (0,6) node[\biasColor] (1) {\gainGraph{-2t}   {s-2t}{s-t}{s}{t}{STUVWXYZ}};
      \draw[shift={(0,3)}, xscale=3]
      (-1.5,0) node             (svz)  {\gainGraph{-2t}   {c-t} {c}  {s}{t}{SVZ}} %
      ( -.5,0) node[\biasColor] (tz)   {\gainGraph{c-s-t} {c-t} {c}  {s}{t}{TZ}}
      (  .5,0) node             (uz)   {\gainGraph{2(c-s)}{c-t} {c}  {s}{t}{UZ}} %
      ( 1.5,0) node             (wxyz) {\gainGraph{a}     {s-2t}{s-t}{s}{t}{WXYZ}}
      ;
      \draw
      (0,0) node                (0)    {\gainGraph{a}     {c-t} {c}  {s}{t}{Z}}
      ;
      \draw
      (1) edge (svz.north) edge (tz.north) edge (uz.north) edge (wxyz.north)
      (0) edge (svz.south) edge (tz.south) edge (uz.south) edge (wxyz.south)
      ;
    \end{tikzpicture}
    \caption[Modular circuits with balance \( \linearClassSymbol_2 \).]{%
      \label{figure:bias-2}%
      \( \linearClassSymbol_2: b=c-t \).
    }
  \end{figure}

  \begin{figure}
    \centering
    \begin{tikzpicture}[every node/.style={inner sep=2pt, scale=\scaleParam, \normalColor}]
      \draw
      (0,6)  node[\biasColor] (1)   {\gainGraph{-2t}  {s-2t}{s-t}{s}{t}{STUVWXYZ}};
      \draw[shift={(0,3)}, xscale=3]
      (-2,0) node             (s)   {\gainGraph{c-s-t}{\tfrac{1}{2}(3c-s-t)}{c}{s}{t}{S}}
      (-1,0) node[\biasColor] (tz)  {\gainGraph{c-s-t}{c-t} {c}  {s}{t}{TZ}}
      ( 0,0) node             (uvx) {\gainGraph{-2t}  {b}   {s-t}{s}{t}{UVX}}
      ( 1,0) node             (w)   {\gainGraph{c-s-t}{2c-s}{c}  {s}{t}{W}}
      ( 2,0) node             (y)   {\gainGraph{c-s-t}{s-2t}{c}  {s}{t}{Y}}
      ;
      \draw
      (0,0)  node             (0)   {\gainGraph{c-s-t}{b}   {c}  {s}{t}{\emptyset}};
      \draw
      (1) edge (s.north) edge (tz.north) edge (uvx.north) edge (w.north) edge (y.north)
      (0) edge (s.south) edge (tz.south) edge (uvx.south) edge (w.south) edge (y.south)
      ;
    \end{tikzpicture}
    \caption[Modular circuits with balance \( \linearClassSymbol_3 \).]{%
      \label{figure:bias-3}%
      \( \linearClassSymbol_3: a=c-s-t \).
    }
  \end{figure}

  \begin{figure}
    \centering
    \begin{tikzpicture}[every node/.style={inner sep=2pt, scale=\scaleParam, \normalColor}]
      \draw
      (0,3) node (1) {\gainGraph{-2t}{s-2t}{d-t}{s}{t}{STUVWXYZ}}
      (0,0) node (0) {\gainGraph{c-s-t}{c-t}{c}{s}{t}{TZ}}
      (1) edge (0)
      ;
    \end{tikzpicture}
    \caption[Modular circuits with balance \( \linearClassSymbol_4 \).]{%
      \label{figure:bias-4}%
      \( \linearClassSymbol_4: a=b-s, b=c-t \).
    }
  \end{figure}
}

\subsection{Low dimension}
In dimensions one and two, the conditions obtained in previous sections can be gotten more concretely.
This section re-derives those conditions by directly manipulating the arrangements and configurations.

\subsubsection{Generic Gains in Dimension One}
  Let \( \mc{Q} \) be a set of \( n \) distinct points in \( \R^1 \).
  Note that hyperplanes in \( \R^1 \) are points.
  Then \( (\mc{Q}, \Gamma, g) \) is gain-generic if and only if no two of the hyperplanes coincide.

  Letting each edge \( e \) of \( \Gamma \) be oriented from \( u = u_e \) to \( v = v_e \), we obtain hyperplanes
  \begin{align*}
    g(e)
    &= d^2(x, u) - d^2(x, v) \\
    &= (x - u)^2 - (x - v)^2 \\
    &= (v - u)(2x - (u + v))
      .
  \end{align*}
  Solving for the unknown yields
  \[
    x = \frac{g(e) + v^2 - u^2}{2(v - u)}
    .
  \]
  Thus each pair \( e, f \) of distinct edges of \( \Gamma \) yields one equation of nongenericity,
  \[
    \frac{g(e) + v_e^2 - u_e^2}{2(v_e - u_e)}
    = \frac{g(f) + v_f^2 - u_f^2}{2(v_f - u_f)}
    ,
  \]
  which we rewrite as
  \[
    (v_e - u_e) \cdot \tfrac{-1}{2} (g(f) + v_f^2 - u_f^2) - (v_f - u_f) \cdot \tfrac{-1}{2} (g(e) + v_e^2 - u_e^2) = 0
    .
  \]
  The left side of this final equation is the determinant of the matrix
  \[
    \systemMatrix{\Gamma}{\mc{Q}}{ef}
    =
    \begin{bmatrix}
      v_e - u_e & -\tfrac{1}{2}(g(e) + v_e^2 - u_e^2) \\
      v_f - u_f & -\tfrac{1}{2}(g(f) + v_f^2 - u_f^2)
    \end{bmatrix}
    .
  \]
  Thus we obtain the non-generic equation \( \det \systemMatrix{\Gamma}{\mc{Q}}{ef} = 0 \), i.e., the equation of the hyperplane \( F_{ef} \), as described in \cref{result:gain-generic-forbidden-equations}---note that this required \( 2 = d+1 \) edges exactly.
  The matroid \( M_\Gamma(\mc{Q}) \) is the uniform matroid of rank \( 1 \) on the edges of \( \Gamma \).
  Moreover, we have one hyperplane of nongenericity \( F_{ef} \) in the edge space of \( \Gamma \) for each pair \( e, f \) of distinct edges of \( \Gamma \), i.e., for each circuit of \( M_\Gamma(\mc{Q}) \).

\subsubsection{Generic Gains in Dimension Two}
  The difference when moving to dimension two is instructive:\ our configuration can now have nontrivial projective dependencies.
  Let \( \mc{Q} \) be a configuration in \( \affineSpace^2 \) and let \( \Gamma \) be a graph with edges \( e \) oriented from \( u = u_e \) to \( v = v_e \).
  Parallelisms are possible in \( M_\Gamma(\mc{Q}) \) arising from \( \affineSpan{q_i, q_j} \parallel \affineSpan{q_k, q_l} \) rather than parallel edges of \( \Gamma \).

  If \( e, f \) is a pair of edges in \( \Gamma \) with the line segments determined by \( e \) and \( f \) parallel in \( \affineSpace^2 \), then each hyperplane (here, line) has equation
  \begin{align*}
    g(e)
    &= d^2(x, u) - d^2(x, v) \\
    &= |x - u|^2 - |x - v|^2 \\
    &= (x - u) \cdot (x - u) - (x - v) \cdot (x - v) \\
    &= (|x|^2 - 2 x \cdot u + |u|^2) - (|x|^2 - 2 x \cdot v + |v|^2) \\
    &= 2 x \cdot (v - u) + (|u|^2 - |v|^2)
      .
  \end{align*}
  In this case we obtain
  \[
    2 x \cdot (v - u)
    = g(e) + |v|^2 - |u|^2
    .
  \]
  Thus noting \( v_f - u_f = c_{e, f}(v_e - u_e) \) for \( c_{e, f} \neq 0 \) we have a non-generic equation
  \[
    g(f) + |v_f|^2 - |u_f|^2
    = 2 x \cdot (v_f - u_f)
    = 2 c_{e, f} x \cdot (v_e - u_e)
    = c_{e, f} (g(e) + |v_e|^2 - |u_e|^2)
    .
  \]

  Every other pair of edges must correspond to intersecting lines, and every triple (generically) corresponds to a non-central subarrangement.
  Let \( \alpha = (x(\alpha), y(\alpha)) \) for all \( \alpha \in \R^2 \).
  Given three distinct edges \( a, b, c \) with no pair of corresponding segments parallel, we again obtain a hyperplane \( H(a, b, c) \) in edge-space with equation
  \[
    \det
    \begin{bmatrix}
      x(v_a - u_a) & y(v_a - u_a) & -\tfrac{1}{2}(g(a) + |v_a|^2 - |u_a|^2) \\
      x(v_b - u_b) & y(v_b - u_b) & -\tfrac{1}{2}(g(b) + |v_b|^2 - |u_b|^2) \\
      x(v_c - u_c) & y(v_c - u_c) & -\tfrac{1}{2}(g(c) + |v_c|^2 - |u_c|^2)
    \end{bmatrix}
    =
    0
    .
  \]
  If \( \{a, b, c\} \) is a coherently oriented triangle in \( \Gamma \), then this simplifies to
  \[
    g(a) + g(b) + g(c) = 0
    .
  \]
  Thus, for \( d = 2 \) the triangles of a gain-generic triple \( (\mc{Q}, \Gamma, g) \) are all unbalanced.

\subsection{Equivalent Representations}
This section uses \cref{result:pythagorean-arrangement-combinatorics} to change the structure of the graph used to represent a given Pythagorean arrangement.
Our main goal is to make the combinatorics of the arrangement reflected in the combinatorics of the graph to the greatest extent possible.
For example, we can use the result to force all parallel hyperplanes to be hyperplanes of parallel edges.
At the same time, our method can be used to move the reference points, provided all the directions of the original arrangement are represented in the new configuration.
As a result, we can simplify the given data and turn geometric structures into combinatorial ones.

Two triples \( (\mc{Q}, \Gamma, g) \) and \( (\mc{Q}', \Gamma', g') \) are \emph{equivalent} when \( \pythag{\mc{Q}}{\Gamma}{g} = \pythag{\mc{Q}'}{\Gamma'}{g'} \).

\begin{rmk}
  There are obvious specializations of this definition to linear equivalence and projective equivalence of the underlying arrangements.
  We leave the study of these specializations as future work.
\end{rmk}

\subsubsection{Gain transportation}
Given \( (\mc{Q}, \Gamma) \) and \( (\mc{Q}', \Gamma') \), where edge \( uv \in E(\Gamma) \) is associated to edge \( u'v' \in \Gamma' \) subject to the constraint \( q_u q_v \parallel q_{u'} q_{v'} \), we may use \( \mc{F}(\mc{Q} \cup \mc{Q}', \Gamma \sqcup \Gamma') \) to transport any gain function \( g \) on \( \Gamma \) to a gain function \( g' \) on \( \Gamma' \) for which \( \mc{H}(\mc{Q}, \Gamma, g) = \mc{H}(\mc{Q}, \Gamma, g) \).
In particular, by assumption the pair \( \{uv, u'v'\} \) is a circuit of \( M_{\Gamma \sqcup \Gamma'}(\mc{Q} \cup \mc{Q'}) \) the equation for \( F_{uv, u'v'} \) determines the gain on \( u'v' \) from that on \( uv \).
We thus obtain the following.

\begin{prop}
  For all \( (\mc{Q}, \Gamma, g) \), there is an equivalent \( (\mc{Q}', \Gamma', g') \) for which \( H_e \parallel H_f \) in \( \pythag{\mc{Q}}{\Gamma}{g} \) implies \( e \parallel f \) in \( \Gamma' \).
  \noProof
\end{prop}
As a simple application, we obtain the following simple classification of Pythagorean arrangements in dimension one.
\begin{cor}
  Every triple \( (\mc{Q}, \Gamma, g) \) is canonically equivalent to a unique triple of the form \( (\{0, 1\}, \Gamma', g') \) for which \( E(\Gamma') = E(\Gamma) \).
  Moreover, the derived arrangement \( \mc{F}(\{0, 1\}, \Gamma') \) in edge space (having fixed an orientation on \( E(\Gamma') \)) is the graphic arrangement of the complete graph with vertices \( E(\Gamma) \).
\end{cor}
\begin{proof}
  Every hyperplane of \( \mc{F}(\{0, 1\}, \Gamma') \) has the form \( g_e - g_f = 0 \).
\end{proof}

Beyond merely telegraphing parallelism in the combinatorics, it is possible to highlight the contribution of any given circuit of \( M_\Gamma(\mc{Q}) \) as a circle of an equivalent triple.
\begin{prop}
  For all \( C \in M_\Gamma(\mc{Q}) \), there are \( \Gamma' \) and \( \mc{Q}' \) for which \( M_{\Gamma'}(\mc{Q}') = M_\Gamma(\mc{Q}) \) and \( C \) is a circle of \( \Gamma' \).
\end{prop}
\begin{proof}
  As \( C \) is a circuit of \( M_\Gamma(\mc{Q}) \), there is a linear dependence relation \( \sum_{uv \in C} r_{uv} (q_v - q_u) = 0 \).
  Choosing an arbitrary point \( q_0' \) and any order of \( C \), define \( q_{i+1}' = q_i' + r_{v_i v_{i+1}} (q_{v_{i+1}} - q_{v_i}) \).
  By construction, \( q_0' = q_n' \) where \( C \) has length \( n \); moreover \( q_{i+1}' - q_i' \) is a vector parallel to \( q_{v_{i+1}} - q_{v_i} \).
  Finally, use the technique of this section to transport edges of \( C \) to the newly created reference points.
\end{proof}

On the other hand, it is also possible to suppress the contribution of all circles.
\begin{prop}
  For all configurations \( \mc{Q} \), graphs \( \Gamma \), and trees \( T \) with \( |E(T)| = |E(\Gamma)| \), there is a set of reference points \( \mc{Q}' \) for which \( M_T(\mc{Q}') = M_\Gamma(\mc{Q}) \).
\end{prop}
\begin{proof}
  Choose an arbitrary root-vertex for \( T \) and a reference point corresponding to this vertex.
  To place a new hyperplane \( H_e \) for which an end, \( q \), has already been placed, construct a new reference point at \( q' = q + c v \) where \( v \) is the normal vector of \( H \), scaling \( v \) by \( c \) if necessary so that \( q' \) does not coincide with any previously placed reference point.
  Finally, having constructed \( \mc{Q}' \) as above, transport the gains of \( \Gamma \) to gains in \( T \) via the technique from this section.
\end{proof}

We shall soon see examples showing the impossibility of confining the combinatorics of a Pythagorean arrangement solely to the combinatorics of the corresponding gain graph.

\subsubsection{Reworking \cref{example:forbidden-gains}}%
  \label{ex:theta-reworked}
  We use the technique described above to encode the parallelism of hyperplanes \( h_{12} \) and \( h_{34} \) of \cref{example:forbidden-gains} into the graph structure in \cref{fig:theta-reworked-graph}.
  Note that \( \Gamma \) and \( \Gamma' \) both have four vertices, the minimum possible given that the original arrangement has four hyperplane directions.
  In \cref{fig:theta-reworked-points} we show the positions of the points in the configuration \( \mc{Q} \cup \mc{Q}' \).
  \begin{figure}
    \centering
    \newcommand{\oldColour}{Green}
    \newcommand{\newColour}{Blue}
    \begin{subfigure}[b]{.4 \textwidth}
      \centering
      \begin{tikzpicture}[scale=.8, transform shape]
        \begin{scope}[every node/.style={circle, fill, inner sep=1.5pt}]
          \draw[\oldColour]
          (0, 0) node[label={270:$1$}] (A) {}
          (3, 0) node[label={270:$2$}] (B) {}
          (3, 3) node[label={ 90:$3$}] (C) {}
          (0, 3) node[label={ 90:$4$}] (D) {}
          ;
          \draw[\newColour]
          (4, 0) node[label={270:$5$}] (E) {}
          (7, 3) node[label={ 90:$6$}] (F) {}
          (7, 0) node[label={270:$7$}] (G) {}
          (4, 3) node[label={ 90:$8$}] (H) {}
          ;
        \end{scope}
        \begin{scope}[
            every node/.style={
              inner sep=1.5pt,
              inner sep=2pt,
              fill=white,
              sloped
            },
          ]
          \draw[\oldColour]
          (A) edge node {$g_{12}$} (B)
          (A) edge node {$g_{13}$} (C)
          (A) edge node {$g_{14}$} (D)
          (B) edge node {$g_{23}$} (C)
          (C) edge node {$g_{34}$} (D)
          ;
          \draw[\newColour]
          (E) edge             node {$g_{14} - 2$}   (F)
          (E) edge[bend left]  node {$-g_{34}$}      (G)
          (E) edge[bend right] node {$g_{12}/2$}     (G)
          (E) edge             node {$-g_{13} - 20$} (H)
          (G) edge             node {$2 - g_{23}$}   (F)
          ;
        \end{scope}
      \end{tikzpicture}
      \subcaption{%
        \label{fig:theta-reworked-graph}%
        \( \textcolor{\oldColour}{\Gamma} \cup \textcolor{\newColour}{\Gamma'} \)
      }
    \end{subfigure}
    \begin{subfigure}[b]{.55 \textwidth}
      \centering
      \begin{tikzpicture}
        \draw[
          point/.style={circle, inner sep=1.5pt, fill=\oldColour},
          bg/.style={inner sep=1ex, \oldColour},
        ]
        (0, 0) node[point] {} node[left , bg] {$q_{1} = (0,0)$}
        (4, 0) node[point] {} node[right, bg] {$q_{2} = (4,0)$}
        (3, 2) node[point] {} node[below, bg] {$q_{3} = (3,2)$}
        (1, 2) node[point] {} node[left , bg] {$q_{4} = (1,2)$}
        ;
        \draw[
          point/.style={circle, inner sep=1.5pt, fill=\newColour},
          bg/.style={inner sep=1ex, \newColour},
        ]
        (1, 0) node[point] {} node[above, bg] {$q_{5} = (1,0)$}
        (2, 2) node[point] {} node[above, bg] {$q_{6} = (2,2)$}
        (3, 0) node[point] {} node[below, bg] {$q_{7} = (3,0)$}
        (4, 2) node[point] {} node[right, bg] {$q_{8} = (4,2)$}
        ;
      \end{tikzpicture}
      \subcaption{%
        \label{fig:theta-reworked-points}%
        \( \textcolor{\oldColour}{\mc{Q}} \cup \textcolor{\newColour}{\mc{Q}'} \)
      }
    \end{subfigure}
    \caption{%
      \label{fig:theta-reworked}%
      Reworking \cref{example:forbidden-gains}.
    }
  \end{figure}
  Using the given coordinates in tandem with \cref{equation:forbidden-hyperplane-definition}, we may convert gains on edges of \( \Gamma \) to gains on edges of \( \Gamma' \).
  The relevant equations are given below, and the results are displayed in \cref{fig:theta-reworked-graph}.
  \begin{align*}
    g_{12} - 2g_{57} & = 0
    &
    g_{13} + g_{58} & = -20
    &
    g_{14} - g_{56} & = 2
    &
    g_{23} + g_{67} & = 2
    &
    g_{34} + g_{57} & = 0
  \end{align*}

  Notice that while the new graph reflects more of the combinatorics of the arrangement, it still does not detect some of the possible triple points via circles in the graph.
  This is expected, and it is not possible to reduce all of the combinatorics to that of the graph; indeed, the matroid \( M_\Gamma(\mc{Q}) = M_{\Gamma'}(\mc{Q}') \) is not graphic as it contains a copy of the four point line, i.e., the uniform matroid \( U_{2, 2} \) of rank two with four elements.
  To see such a four point line, consider any four edges of \( \Gamma' \), only one of which is selected from the digon.
  Note that there are (by construction) four circuits of \( M_{\Gamma'}(\mc{Q}') \) with size three on these four edges, corresponding with the four possible ways of making a single triple-point in the restricted arrangement.
  Hence \( M_{\Gamma'}(\mc{Q}')|S = U_{2, 2} \).

\subsection{The Pappos Line Arrangement}
\label{ex:pappos-example}
{
  \newcommand{\pointConfiguration}[1]{
    \foreach \l/\a/\b/\ang in {
      1/4/0/below,
      2/0/1/below,
      3/2/1/above,
      4/3/1/below,
      5/5/1/above,
      6/1/2/left,
      7/1/3/above,
      8/4/3/above
    }{#1}
  }
  \newcommand{\graphEdges}[1]{
    \foreach \a/\b/\loD/\hiD in {
      1/4/4/6,
      1/5/2/8,
      2/3/3.5/3.5,
      2/6/8/2,
      2/8/2.25/1.25,
      3/6/6/4,
      4/5/3.5/3.5,
      5/7/1.25/2.25,
      7/8/3.675/1}{#1}
  }
  Consider the point configuration and graph of \cref{figure:pappos-arrangement-obstruction}.
  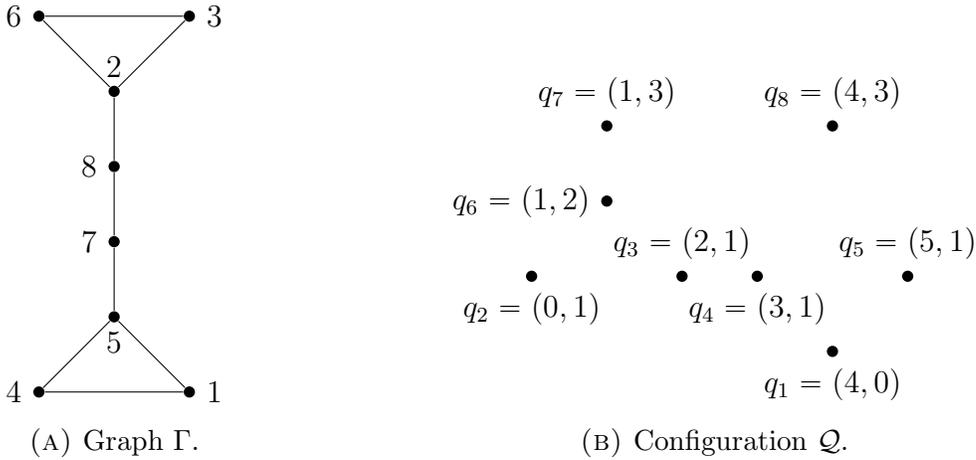
\begin{figure}
    \centering
    \begin{subfigure}[b]{.3\textwidth}
      \centering
      \begin{tikzpicture}[scale=1, transform shape]
        \draw[every node/.style={inner sep=1.5pt, fill, circle}]
        ( 1,-2.5) node[label={right:$1$}] (1) {}
        ( 0, 1.5) node[label={above:$2$}] (2) {}
        ( 1, 2.5) node[label={right:$3$}] (3) {}
        (-1,-2.5) node[label={left :$4$}] (4) {}
        ( 0,-1.5) node[label={below:$5$}] (5) {}
        (-1, 2.5) node[label={left :$6$}] (6) {}
        ( 0,-0.5) node[label={left :$7$}] (7) {}
        ( 0, 0.5) node[label={left :$8$}] (8) {}
        ;
        \graphEdges{\draw (\a) -- (\b);}
      \end{tikzpicture}
      \subcaption{Graph \( \Gamma \).}
    \end{subfigure}
    \begin{subfigure}[b]{.65\textwidth}
      \centering
      \begin{tikzpicture}[scale=1, transform shape]
        \pointConfiguration{
          \coordinate (\l) at (\a,\b);
          \node[circle, fill, inner sep=1.5pt, label={\ang:$q_\l=(\a,\b)$}] at (\l) {};
        }
      \end{tikzpicture}
      \subcaption{Configuration \( \mc{Q} \).}
    \end{subfigure}
    \caption[Pappos line arrangement as a Pythagorean arrangement.]{%
      \label{figure:pappos-arrangement-obstruction}%
      An obstruction to immediate combinatorialization.
    }
  \end{figure}
  Then \( \pythag{\mc{Q}}{\Gamma}{\mathbf{0}} \) is the Pappos line arrangement of \cref{figure:pappos-line-arrangement}.
  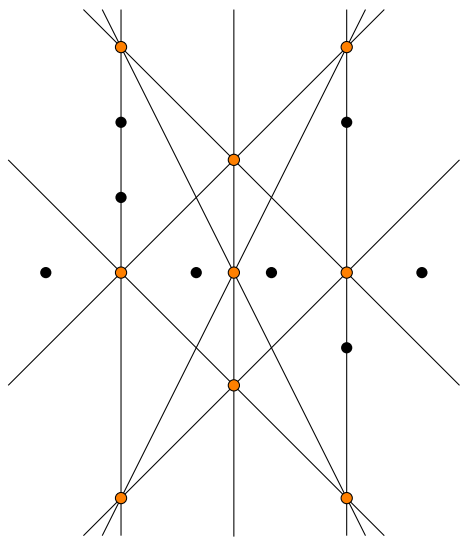
\begin{figure}
    \centering
    \begin{tikzpicture}[scale=1]
      \pointConfiguration{
        \coordinate (\l) at (\a,\b);
        \node[circle, fill, inner sep=1.5pt] at (\l) {};
      }
      \graphEdges{
        \path
        (\a) -- (\b) coordinate[midway] (\a\b)
        ($(\a\b)!\loD!90:(\a)$)  coordinate (lo\a\b)
        ($(\a\b)!\hiD!270:(\a)$) coordinate (hi\a\b)
        (lo\a\b) edge (hi\a\b)
        ;
      }
      \draw[every node/.style={circle, fill, inner sep=1.5pt, orange, draw=black}]
      (intersection of lo14--hi14 and lo15--hi15) node {}
      (intersection of lo14--hi14 and lo23--hi23) node {}
      (intersection of lo14--hi14 and lo78--hi78) node {}
      (intersection of lo15--hi15 and lo23--hi23) node {}
      (intersection of lo15--hi15 and lo78--hi78) node {}
      (intersection of lo23--hi23 and lo26--hi26) node {}
      (intersection of lo26--hi26 and lo45--hi45) node {}
      (intersection of lo28--hi28 and lo78--hi78) node {}
      (intersection of lo45--hi45 and lo57--hi57) node {}
      ;
    \end{tikzpicture}
    \caption[Pappos line arrangement.]{%
      \label{figure:pappos-line-arrangement}
      Pappos line arrangement.%
    }
  \end{figure}
  We view this as an obstruction to immediate combinatorial abstraction.
  In particular, without the geometric data encoded in \( M_\Gamma(\mc{Q}) \) we would not detect the three broken triple-points obtained by wiggling the gain on any one edge, as there is no structure in \( \Gamma \) alone that detects this quirk.

  However, this issue is primarily caused by the fact that all of the parallelism present in the arrangement is caused by geometric parallelism, which is not reflected in the graph.
  In particular, there are fundamentally only five directions used from among the reference points, so we ought be able to significantly reduce the number of reference points and at the same time force the combinatorics of the arrangement to be better reflected by that of the new graph.
  We exhibit an equivalent triple \( (\mc{Q}', \Gamma', g') \) in \cref{figure:pappos-arrangement-reworked}.
  \begin{figure}
    \centering
    \newcommand{\xA}{-5/2}
    \newcommand{\yA}{0}
    \newcommand{\xB}{5/2}
    \newcommand{\yB}{0}
    \newcommand{\xC}{0}
    \newcommand{\yC}{5/4}
    \newcommand{\xD}{0}
    \newcommand{\yD}{-5/2}
    \newcommand{\gABa}{0}
    \newcommand{\gABb}{-15}
    \newcommand{\gABc}{15}
    \newcommand{\gACa}{75/16}
    \newcommand{\gADa}{15/2}
    \newcommand{\gADb}{-15/2}
    \newcommand{\gBCa}{75/16}
    \newcommand{\gBDa}{-15/2}
    \newcommand{\gBDb}{15/2}
    \pgfmathsetmacro{\xMin}{-2.5}
    \pgfmathsetmacro{\xMax}{2.5}
    \pgfmathsetmacro{\yMin}{-4}
    \pgfmathsetmacro{\yMax}{4}
    \newcommand{\cABa}{orange}
    \newcommand{\cABb}{blue}
    \newcommand{\cABc}{green}
    \newcommand{\cACa}{orange}
    \newcommand{\cADa}{green}
    \newcommand{\cADb}{blue}
    \newcommand{\cBCa}{orange}
    \newcommand{\cBDa}{green}
    \newcommand{\cBDb}{blue}
    \begin{subfigure}[b]{.35 \textwidth}
      \centering
      \begin{tikzpicture}[scale=1, transform shape]
        \draw[every node/.style={circle, fill, inner sep=1.5pt}]
        (-2, 0) node[label={180:$1$}] (A) {}
        ( 2, 0) node[label={  0:$2$}] (B) {}
        ( 0, 2) node[label={ 90:$3$}] (C) {}
        ( 0,-2) node[label={270:$4$}] (D) {}
        ;
        \draw[
          every node/.style={inner sep=1.5pt},
          lft1/.style={bend left=4ex},
          rht1/.style={bend right=4ex},
          lft2/.style={bend left=10ex},
          rht2/.style={bend right=10ex},
        ]
        (A) edge       node[fill=white]{$\gABa$} node {$\gABa$} (B)
        (A) edge[lft1] node[fill=white]{$\gABb$} node {$\gABb$} (B)
        (A) edge[rht1] node[fill=white]{$\gABc$} node {$\gABc$} (B)
        (A) edge       node[fill=white]{$\gACa$} node {$\gACa$} (C)
        (A) edge       node[fill=white]{$\gADa$} node {$\gADa$} (D)
        (A) edge[rht2] node[fill=white]{$\gADb$} node {$\gADb$} (D)
        (B) edge       node[fill=white]{$\gBCa$} node {$\gBCa$} (C)
        (B) edge[lft2] node[fill=white]{$\gBDa$} node {$\gBDa$} (D)
        (B) edge       node[fill=white]{$\gBDb$} node {$\gBDb$} (D)
        ;
      \end{tikzpicture}
      \subcaption{Simplified graph \( \Gamma' \).}
    \end{subfigure}
    \begin{subfigure}[b]{.55 \textwidth}
      \centering
      \begin{tikzpicture}[scale=1, transform shape]
        \draw[every node/.style={fill, inner sep=1.5pt}]
        (\xA, \yA) node[circle, label={270:$q_1' = (\xA, \yA)$}] (A) {}
        (\xB, \yB) node[circle, label={270:$q_2' = (\xB, \yB)$}] (B) {}
        (\xC, \yC) node[circle, label={ 90:$q_3' = (\xC, \yC)$}] (C) {}
        (\xD, \yD) node[circle, label={270:$q_4' = (\xD, \yD)$}] (D) {}
        ;
      \end{tikzpicture}
      \subcaption{Simplified configuration \( \mc{Q}' \).}
    \end{subfigure}
    \caption[Pappos line arrangement as a better Pythagorean arrangement.]{%
      \label{figure:pappos-arrangement-reworked}%
      The modified gain graph and point configuration.
    }
  \end{figure}
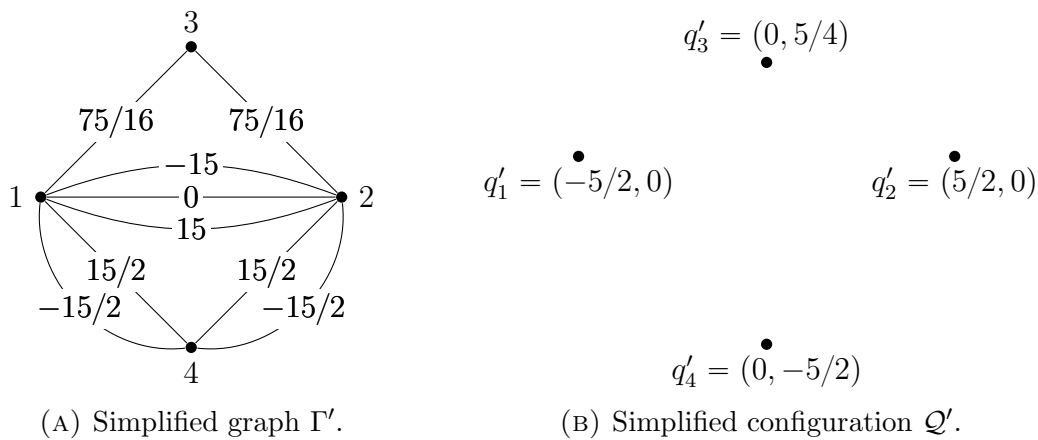
  As with \cref{ex:theta-reworked}, this new representation uses the minimum possible number of points.
  Its matroid at infinity also includes a four-point line, and thus the graph structure cannot capture the combinatorics of the resulting arrangement without the geometric data provided by \( M_\Gamma(\mc{Q}) \).
}

\printbibliography
\end{document}